\documentclass[12pt]{amsart}
\usepackage{colordvi}
\usepackage{graphicx,xspace}
\usepackage{bm}
\usepackage{empheq,bbm}
\usepackage{amssymb}
\usepackage[only,llbracket,rrbracket,bbslash]{stmaryrd}
\usepackage{pstricks,pst-node}
\usepackage[width=.90\textwidth,font=small,labelfont=sc]{caption}
\usepackage[bookmarks=false,colorlinks,linkcolor=mylinkcolor,citecolor=mycitecolor]{hyperref} 
\usepackage{bookmark}
	\definecolor{mycitecolor}{rgb}{.2,.7,.2}
	\definecolor{mylinkcolor}{rgb}{0,0,.8}

\newfam\cyrfam
\font\cyr=wncyr10 at 8pt   
\newfam\cyrfam
\font\Cyr=wncyr10 at 12pt   

\newcommand{\ch}{\mbox{\cyr Ch}}
\newcommand{\CH}{\mbox{\Cyr Ch}}

\numberwithin{equation}{section}
\newtheorem{theorem}{Theorem}[section]
\newtheorem{prop}[theorem]{Proposition}
\newtheorem{lemma}[theorem]{Lemma}
\newtheorem{corollary}[theorem]{Corollary}

\newtheorem{definition}[theorem]{Definition}
\newtheorem{ex}[theorem]{Example}
\newtheorem{rem}[theorem]{Remark}

\newenvironment{example}[1][]{\rm\begin{ex}[#1]\rm}{\end{ex}}
\newenvironment{remark}[1][]{\rm\begin{rem}[#1]\rm}{\end{rem}}

\newpsobject{showgrid}{psgrid}{%
  subgriddiv=1,griddots=10,gridlabels=6pt}
\newif\ifhrule\hrulefalse

\newcounter{FNC}[page]
\def\fauxfootnote#1{{\addtocounter{FNC}{2}$^\fnsymbol{FNC}$%
     \let\thefootnote\relax\footnotetext{$^\fnsymbol{FNC}$\Magenta{#1}}}}

\PII{}
\copyrightinfo{}{}

\def\Q{\mathbb Q}

\newcommand{\chains}{\ch}
\newcommand{\Chains}{\CH}

\def\s{\mathfrak S}
\def\m{\mathcal M}
\def\y{\mathcal Y}

\def\calB{{\mathcal B}}
\def\calS{{\mathcal S}}

\def\ssym{\mathfrak SSym}
\def\msym{\mathcal MSym}
\def\ysym{\mathcal YSym}

\newcommand{\setS}{{\sf S}}
\newcommand{\setT}{{\sf T}}

\def\id{\mathbbm 1}
\def\btau{\bm\tau}
\def\bbeta{\bm\beta}
\def\bphi{\bm\phi}
\def\brho{\bm\rho}

\newcommand{\Placeholder}%
	{\raisebox{.09ex}{\footnotesize $\bullet$}}
\newcommand{\placeholder}%
	{\raisebox{.04ex}{\tiny $\bullet$}}
	
\def\bb{\boldbullet}
	\def\boldbullet{{\!\mbox{\LARGE$\mathbf\cdot$}}}  

\def\psplit{\stackrel{\curlyvee}{\to}}
\def\rsplit{\stackrel{\curlyvee_{\!\!+}}{\longrightarrow}}

\newcommand{\Max}[1][]{\mathsf{max}\if#1\else{\left(#1\right)}\fi}
\newcommand{\Min}[1][]{\mathsf{min}\if#1\else{\left(#1\right)}\fi}
\newcommand{\bMax}[1][]{{\mathsf{max}}\if#1\else{\left(#1\right)}\fi}

\def\onto{\twoheadrightarrow}

\newcommand{\enum}[2]{\mathbf{#2}({#1})}
\newcommand{\demph}[1]{{\Blue{{\sl #1}}}}

\author{Stefan Forcey} \address[S. Forcey]{
	Department of Mathematics\\
         Tennessee State University\\
          3500 John A Merritt Blvd\\
         Nashville, Tennessee \ 37209  
	} 
	\email{sforcey@tnstate.edu}  \urladdr{http://faculty.tnstate.edu/sforcey/}

\author{Aaron Lauve} \address[A. Lauve]{
	Department of Mathematics\\ 
         Texas A\&M University\\
         MS 3368\\
         College Station, Texas \ 77843 
         }
	\email{lauve@math.tamu.edu}  \urladdr{http://www.math.tamu.edu/\~{}lauve}

\author{Frank Sottile} \address[F. Sottile]{
	Department of Mathematics\\ 
         Texas A\&M University\\
         MS 3368\\
         College Station, Texas \ 77843 
         }
	\thanks{Sottile supported by NSF grant DMS-0701050}  
	\email{sottile@math.tamu.edu}  \urladdr{http://www.math.tamu.edu/\~{}sottile}

\title[Hopf structures on the multiplihedra]{Hopf structures on the multiplihedra} 

\keywords{multiplihedron, permutations, permutahedron, associahedron, binary trees, Hopf algebras}

\begin{document} 

\begin{abstract} 
 We investigate algebraic structures that can be placed on vertices of the
 multiplihedra, a family of polytopes originating in the study of higher
 categories and homotopy theory. 
 Most compelling among these are two distinct structures of a Hopf module 
 over the Loday--Ronco Hopf algebra. 
\end{abstract}

\maketitle



\section*{Introduction}\label{sec: intro}

The permutahedra $\s_\bb$ form a family of highly symmetric polytopes that have been of
interest since their introduction by Schoute in 1911~\cite{Schoute:1911}.
The associahedra $\y_\bb$ are another family of polytopes that were introduced by Stasheff as
cell complexes in 1963~\cite{Sta:1963}, and with the permutahedra were studied from the
perspective of monoidal categories and $H$-spaces~\cite{Mil:1966} in the 1960s.
Only later were associahedra shown to be polytopes~\cite{haiman,lee,milnor}.
Interest in these objects was heightened in the 1990s, when Hopf algebra structures  
were placed on them in work of Malvenuto, Reutenauer, Loday, Ronco, Chapoton, and
others~\cite{Cha:2000,LodRon:1998,MalReu:1995}.  
More recently, the associahedra were shown to arise in Lie theory through work of Fomin and
Zelevinsky on cluster algebras~\cite{FomZel:2003}.

We investigate Hopf structures on another family of polyhedra, the multiplihedra, $\m_\bb$.
Stasheff introduced them in the context of maps preserving higher homotopy
associativity~\cite{Sta:1970} and described their 1-skeleta. 
Boardman and Vogt~\cite{BoaVog:1973}, and then Iwase and Mimura~\cite{IwaMim:1989} described
the multiplihedra as cell complexes, and only recently were they shown to be convex
polytopes~\cite{For:2008}. 
These three families of polytopes are closely related.
For each integer $n\geq 1$, the permutahedron $\s_n$, multiplihedron $\m_n$, and associahedron
$\y_n$ are polytopes of dimension $n{-}1$ with natural cellular surjections
$\s_n\twoheadrightarrow\m_n\twoheadrightarrow\y_n$, which we illustrate when $n=4$.
 \[
  \raisebox{-9ex}{\includegraphics[height=3.2cm]{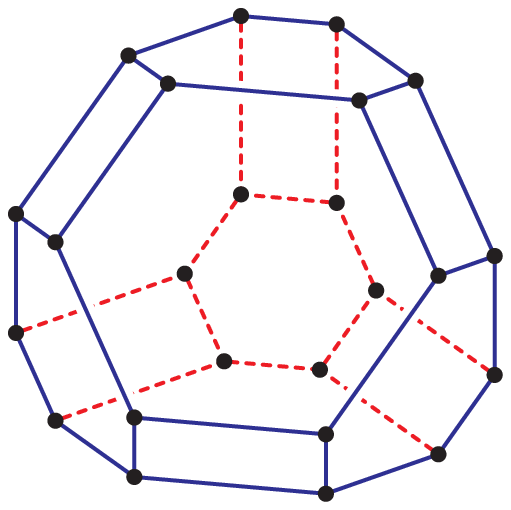}
    \raisebox{9ex}{\  \ $\relbar\joinrel\twoheadrightarrow$\  \ \ }
   \includegraphics[height=3.2cm]{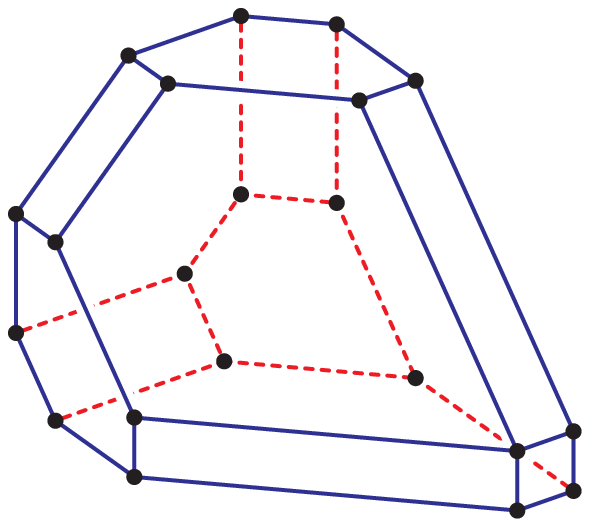}
     \raisebox{9ex}{\ \ $\relbar\joinrel\twoheadrightarrow$\ \ \ }
   \includegraphics[height=3.2cm]{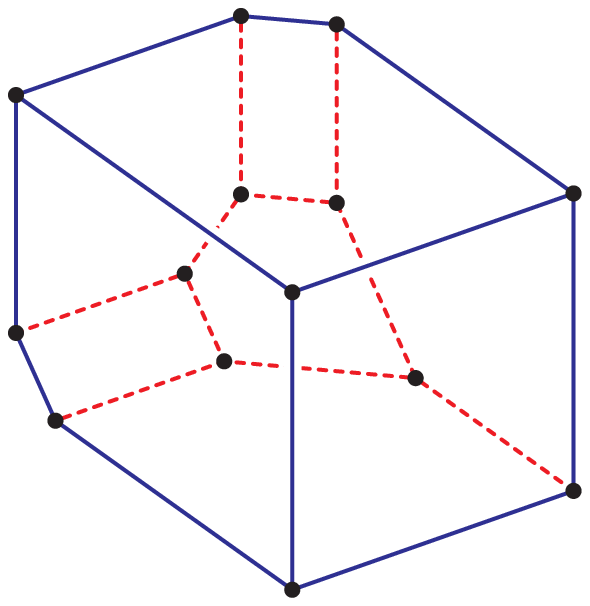}}
 \]
%
%

The faces of these polytopes are represented by different flavors of planar trees;
permutahedra by ordered trees (set compositions), multiplihedra by bi-leveled trees 
(Section~\ref{sec: bi-leveled trees}), and associahedra by planar trees. 
The maps between them forget the additional structure on the trees. 
These maps induce surjective maps of graded vector spaces spanned by
the vertices, which are binary trees.
The span $\ssym$ of ordered trees forms the Malvenuto-Reutenauer Hopf
algebra~\cite{MalReu:1995} and the span $\ysym$ of planar binary trees forms the   
Loday-Ronco Hopf algebra~\cite{LodRon:1998}.
The algebraic structures of multiplication and comultiplication on $\ssym$ and $\ysym$ are
described in terms of geometric operations on trees
and the composed surjection $\tau\colon\s_\bb\twoheadrightarrow\y_\bb$ gives a surjective
morphism $\btau\colon\ssym\twoheadrightarrow\ysym$ of Hopf algebras.

We define $\msym$ to be the vector space spanned by the
vertices of all multiplihedra.
The factorization of $\tau$ induced by the maps of polytopes,
$\ssym\twoheadrightarrow\msym\twoheadrightarrow\ysym$, does not endow $\msym$ with the
structure of a Hopf algebra.
Nevertheless, some algebraic structure does survive the factorization. 
We show in Section
\ref{sec: msym} that $\msym$ is an algebra, which is simultaneously a $\ssym$-module and a
$\ysym$--Hopf module algebra, and the maps preserve these structures.

We perform a change of basis in $\msym$ using M\"obius inversion that illuminates 
its comodule structure.  
Such changes of basis helped to understand the coalgebra structure of
$\ssym$~\cite{AguSot:2005} and of $\ysym$~\cite{AguSot:2006}.
Section~\ref{sec: variations}  discusses a second $\ysym$ Hopf module structure that may 
be placed on the positive part $\msym_+$ of $\msym$. 
This structure also arises from polytope maps between $\s_\bb$ and $\y_\bb$, but not directly 
from the algebra structure of $\ssym.$ 
M\"obius inversion again reveals an explicit basis of $\ysym$ coinvariants 
in this alternate setting.

\section{Basic Combinatorial Data}\label{sec: basic}

The structures of the Malvenuto-Reutenauer and Loday-Ronco algebras are related to the weak
order on ordered trees and the Tamari order on planar trees.
There are natural maps between the weak and Tamari orders which induce a 
morphism of Hopf algebras.
We first recall these partial orders and then the basic structure of these Hopf algebras.
In Section~\ref{sec: poset lemma} we establish a formula involving the M\"obius functions of
two posets related by an interval retract.
This is a strictly weaker notion than that of a Galois correspondence, which was used to study
the structure of the Loday-Ronco Hopf algebra.

\subsection{$\s_\bb$ and $\y_\bb$}\label{sec: s and y} 
The 1-skeleta of the families of polytopes $\s_\bb,\m_\bb$, and $\y_\bb$ are Hasse diagrams
of posets whose structures are intertwined with the algebra structures we study.
We use the same notation for a polytope and its poset of vertices. 
Similarly, we use the same notation for a cellular surjection of polytopes and the poset map formed by restricting that surjection to vertices. 

\begin{figure}[htb]{\small
\[
  \begin{picture}(214,165)(-20,2)
    \put(0,0){\includegraphics[height=160pt]{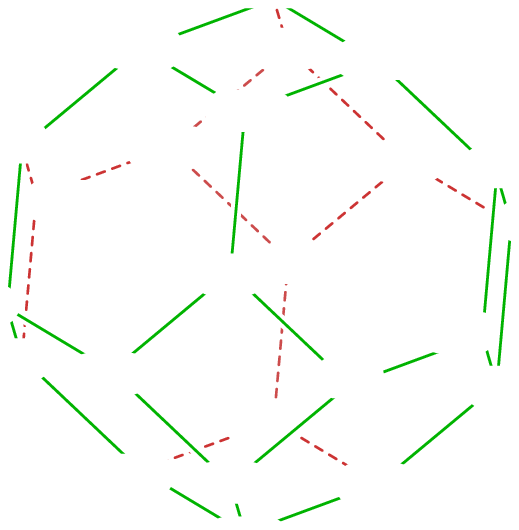}}
                   \put(70,160){\Blue{4321}}
      \put(26,143){\Blue{4312}} \put(72,143){\Blue{4231}}  \put(101,140){\Blue{3421}}
                   \put(62,125){\Blue{3412}}
    \put(-11,114){\Blue{4213}}\put(36,113){\Blue{4132}}\put(109,110){\Blue{3241}}\put(138,107){\Blue{2431}}
    \put(6,98){\Blue{4123}}   \put(149,91){\Blue{2341}}
    \put(-20,66){\Blue{3214}}\put(54,75){\Blue{2413}}\put(79,79){\Blue{3142}}\put(123,58){\Blue{1432}}
    \put(-7,50){\Blue{3124}}\put(22,46){\Blue{2314}}\put(58,31){\Blue{2143}}
          \put(90.5,43.5){\Blue{1423}}\put(135,42){\Blue{1342}}
    \put(28, 16){\Blue{2134}}  \put(56,13){\Blue{1324}}  \put(99,13){\Blue{1243}}
            \put(60,-3){\Blue{1234}}
  \end{picture}
   \qquad
  \begin{picture}(160,165)(-8,8)
   \put(0,6){\includegraphics[height=160pt]{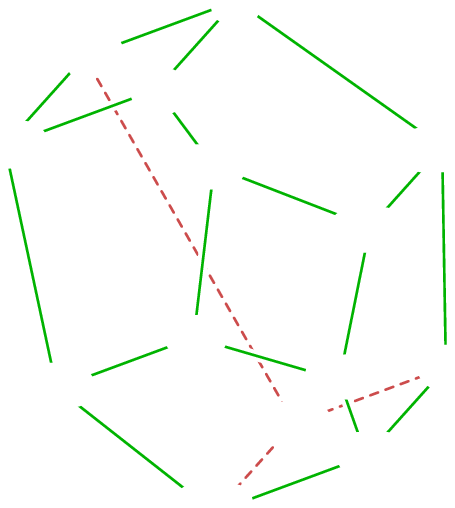}}

   \put( 64,160){\includegraphics[height=12pt]{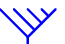}}

   \put( 15,141.5){\includegraphics[height=12pt]{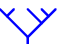}}
   \put( 39,131){\includegraphics[height=12pt]{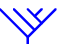}}

   \put(-8.5,115.5){\includegraphics[height=12pt]{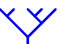}}
   \put( 57.5,109){\includegraphics[height=12pt]{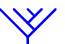}}
   \put(126.5,113){\includegraphics[height=12pt]{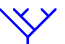}}

   \put(106.5, 88){\includegraphics[height=12pt]{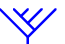}}

   \put(6.5, 39.5){\includegraphics[height=12pt]{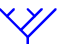}}
   \put(51, 54.5){\includegraphics[height=12pt]{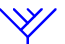}}
   \put( 96, 42){\includegraphics[height=12pt]{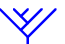}}
   \put(129, 45){\includegraphics[height=12pt]{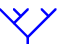}}

   \put( 81, 27.5){\includegraphics[height=12pt]{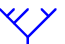}}
   \put(105, 17){\includegraphics[height=12pt]{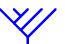}}

   \put( 58,  1.5){\includegraphics[height=12pt]{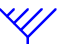}}
  \end{picture}
\]}
 \caption{Weak order on $\s_4$ and Tamari order on $\y_4$}
 \label{Fig:S4Y4}
\end{figure}

For the permutahedron $\s_n$, the corresponding poset is the (left) \demph{weak order}, which we describe in terms
of permutations. 
A cover in the weak order has the form $w\lessdot (k,k{+}1) w$, where $k$ preceeds $k{+}1$
among the values of $w$. 
Figure~\ref{Fig:S4Y4} displays the weak order on $\s_4$. 
We let $\s_0=\{\emptyset\}$, where $\emptyset$ is the empty permutation of $\emptyset$. 

Let \Blue{$\y_n$} be the set of rooted, planar binary trees with $n$ nodes. 
The cover relations in the \demph{Tamari order} on $\y_n$ are obtained by
moving a child node directly above a given node from the left to the right branch above the
given node. 
Thus
 \[
  \raisebox{-7pt}{\includegraphics{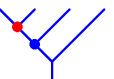}}
      {\color{red} \ \longrightarrow\  }
  \raisebox{-7pt}{\includegraphics{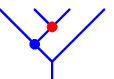}}
    {\color{blue}   \ \longrightarrow\  }
  \raisebox{-7pt}{\includegraphics{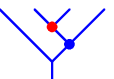}}
      {\color{red} \ \longrightarrow\   }
  \raisebox{-7pt}{\includegraphics{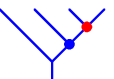}}
 \]
is an increasing chain in $\y_3$ (the moving vertices are marked with  dots).
Figure~\ref{Fig:S4Y4} shows the Tamari order on  $\y_4$. 

The unique tree in $\y_1$ is \includegraphics{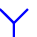}.
Given trees $t_\ell$ and $t_r$, form the tree $t_\ell\vee t_r$ by grafting the root of 
$t_\ell$ (respectively of $t_r$) to the left (respectively right) leaf
of  \includegraphics{figures/1.eps}.
Form the tree $t_\ell\backslash t_r$ by grafting the root of $t_r$ to the rightmost leaf of
$t_\ell$.
For example,
\[
 \begin{picture}(257,45)(-1,0)
  \put( 4,0){$t_\ell$}  \put( 0,12){\includegraphics{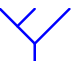}}
  \put(39,0){$t_r$}     \put(35,12){\includegraphics{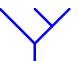}}

  \put(100,0){$t_\ell\vee t_r$}  \put(90,12){\includegraphics{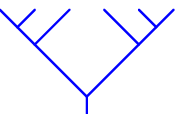}}
 
  \put(160,12){\includegraphics{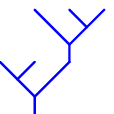}}
  \put(190,17){$=$} 
  \put(205,12){\includegraphics{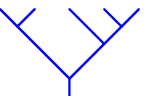}}
  \put(185,0){$t_\ell\backslash t_r$} 
  \put(250,17){\!\!.}
 \end{picture}
\]
Decompositions $t=t_1\backslash t_2$ correspond to pruning $t$ along the right branches from the root. 
A tree $t$ is \demph{indecomposable} if it has no nontrivial decomposition 
$t=t_1\backslash t_2$ with $t_1,t_2\neq\includegraphics{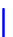}$.
Equivalently, if the root node is the rightmost node of $t$. 
Any tree $t$ is uniquely decomposed $t=t_1\backslash\dotsb\backslash t_m$ into indecomposable
trees $t_1,\dotsc,t_m$.

We define a poset map $\tau\colon \s_n\to\y_n$.
First, given distinct integers $a_1,\dotsc,a_k$, let $\overline{a}\in\s_k$ be
the unique permutation such that $\overline{a}(i)<\overline{a}(j)$ if and only if $a_i<a_j$.
Thus $\overline{4726}=2413$.
Since $\s_0$, $\y_0$, $\s_1$, and $\y_1$ are singletons, we must have
 \begin{eqnarray*}
  \tau\ \colon\ \s_0\ \longrightarrow\ \y_0&\mbox{with}& \tau\ \colon\ 
       \emptyset\ \longmapsto\ \includegraphics{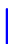}\,,\quad\mbox{and}\\
  \tau\ \colon\ \s_1\ \longrightarrow\ \y_1&\mbox{with}& \tau\ \colon\ 
       1\ \longmapsto\ \includegraphics{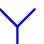}\,.
 \end{eqnarray*}
Let $n>0$ and assume that $\tau$ has been defined on $\s_k$ for $k<n$.
For $w\in\s_n$ suppose that $w(j)=n$, and define
\[
   \tau(w)\ :=\ \tau\bigl(\,\overline{w(1),\dotsc,w(j{-}1)}\,\bigr)\ \vee\ 
                \tau\bigl(\,\overline{w(j{+}1),\dotsc,w(n}\,\bigr)\,.
\]
For example,
\begin{gather*}
   \tau(12)\ =\ \includegraphics{figures/1.d.eps}\
                 \vee\ \includegraphics{figures/0.d.eps} \ =\ 
                \includegraphics{figures/12.d.eps},
   \quad
   \tau(21)\ =\ \includegraphics{figures/0.d.eps}\
                 \vee\ \includegraphics{figures/1.d.eps} \ =\ 
                \includegraphics{figures/21.d.eps},
    \quad\mbox{and}\vspace{-4pt} \\
%
%
   \tau(3421)\ =\ \tau(\overline{3})\ \vee\ \tau(\overline{21})\ =\ 
    \tau(1)\ \vee\ \tau(21)\ =\ 
               \includegraphics{figures/1.d.eps}\
                 \vee\ \includegraphics{figures/21.d.eps} \ =\ 
                \includegraphics{figures/3421.d.eps} .
\end{gather*}

Loday and Ronco~\cite{LodRon:2002} show that the fibers $\tau^{-1}(t)$ of
$\tau$ are intervals in the weak order.
This gives two canonical sections of $\tau$.
For $t\in\y_n$,
 \[
   \Blue{\Min[t]}\ :=\ \min \left\{ w \mid \tau(w) = t \right\} \quad\hbox{ and }\quad 
   \Blue{\Max[t]}\ :=\ \max \left\{ w \mid \tau(w) = t \right\},
 \]
the minimum and maximum in the weak order. 
Equivalently, $\Min[t]$ is the unique $231$-avoiding permutation in $\tau^{-1}(t)$ and
$\Max[t]$ is the unique $132$-avoiding permutation. 
These maps are order-preserving.

The 1-skeleta of $\s_n$ and $\y_n$ form the Hasse diagrams of the weak and Tamari orders,
respectively.
Since $\tau$ is an order-preserving surjection, it induces a cellular map between the 1-skeleta
of these polytopes.
Tonks~\cite{Ton:1997} extended $\tau$ to the faces of $\s_n$, giving a cellular surjection.

The nodes and internal edges of a tree are the Hasse diagram of a poset
with the root node maximal.
Labeling the nodes (equivalently, the gaps between the leaves) of $\tau(w)$ with the 
values of the permutation $w$ gives a linear extension of the node poset of $\tau(w)$,
and all linear
extensions of a tree $t$ arise in this way for a unique permutation in $\tau^{-1}(t)$.
Such a linear extension $w$ of a tree is an \demph{ordered tree} and $\tau(w)$ is the
corresponding unordered tree.
In this way, $\s_n$ is identified with the set of ordered trees with $n$ nodes. 
Here are some ordered trees,
\newcommand{\sn}[1]{{\small\hspace{2pt}$#1$}}
\[
   \begin{picture}(42,33)
   \put(0,0){\includegraphics{figures/3421.d.eps}}
   \put(0,26){\sn{3}} \put(10,26){\sn{4}} \put(20,26){\sn{2}} \put(30,26){\sn{1}}
  \end{picture}
   \qquad
  \begin{picture}(42,33)
   \put(0,0){\includegraphics{figures/3421.d.eps}}
   \put(0,26){\sn{1}} \put(10,26){\sn{4}} \put(20,26){\sn{3}} \put(30,26){\sn{2}}
  \end{picture}
   \qquad
  \begin{picture}(53,38)
   \put(0,0){\includegraphics{figures/34521.d.eps}}
   \put(0,31){\sn{2}} \put(10,31){\sn{3}} \put(20,31){\sn{5}} 
   \put(30,31){\sn{4}}\put(40,31){\sn{1}}
  \end{picture}
   \qquad
  \begin{picture}(52,38)
   \put(0,0){\includegraphics{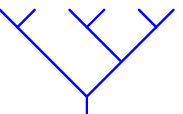}}
   \put(0,31){\sn{2}} \put(10,31){\sn{5}} \put(20,31){\sn{1}} 
   \put(30,31){\sn{4}}\put(40,31){\sn{3}}
  \put(52,1){.}
  \end{picture}
\]

Given ordered trees $u,v$, form the ordered tree $u\backslash v$ by grafting the root of $v$
to the rightmost leaf of $u$,  where the nodes of $u$ are greater than the nodes of $v$, but
the relative orders within $u$ and $v$ are maintained.
Thus we may decompose an ordered tree $w=u\backslash v$ whenever $\tau(w)=r\backslash s$ 
with $\tau(u)=r$, $\tau(v)=s$, and the nodes of $r$ in $w$ precede the nodes of $s$ in $w$.
An ordered tree $w$ is \demph{indecomposable} if it has no nontrivial such decompositions.
Here are ordered trees $u$, $v$ and $u\backslash v$,
\[
  \begin{picture}(42,33)
   \put(0,0){\includegraphics{figures/3421.d.eps}}
   \put(0,26){\sn{1}} \put(10,26){\sn{4}} \put(20,26){\sn{3}} \put(30,26){\sn{2}}
  \end{picture}
   \qquad
  \begin{picture}(42,33)
   \put(0,0){\includegraphics{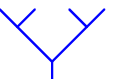}}
   \put(0,21){\sn{1}} \put(10,21){\sn{3}} \put(20,21){\sn{2}}
  \end{picture}
   \qquad\quad
  \begin{picture}(57,53)
   \put(0,0){\includegraphics{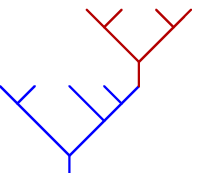}}
   \put( 0,26){\sn{4}} \put(10.5,26){\sn{7}} \put(20,26){\sn{6}}\put(30,26){\sn{5}}
   \put(25,48){\sn{1}} \put(35,48){\sn{3}} \put(45,48){\sn{2}}
  \end{picture}
   \quad\raisebox{15pt}{$=$}\quad
  \begin{picture}(72,48)
   \put( 0,0){\includegraphics{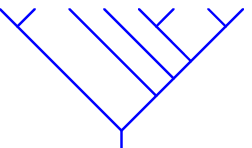}}
   \put( 0,41){\sn{4}} \put(10,41){\sn{7}} \put(20,41){\sn{6}}\put(30,41){\sn{5}}
   \put(40,41){\sn{1}} \put(50,41){\sn{3}} \put(60,41){\sn{2}}
   \put(70,1){.}
  \end{picture}
\]

We may \demph{split} an ordered tree $w$ along a leaf to obtain either an
ordered forest (where the nodes in the forest are totally ordered) or a pair of ordered trees,
\[
   \begin{picture}(52,49)
   \put(0,0){\includegraphics{figures/45231.d.eps}}
   \put(0,31){\sn{2}} \put(10,31){\sn{5}} \put(20,31){\sn{1}} 
   \put(30,31){\sn{4}}\put(40,31){\sn{3}}
   \put(30,49){\vector(0,-1){15}}
  \end{picture}
   \ \raisebox{12pt}{\large$\leadsto$}\ 
  \begin{picture}(52,49)
   \put(0,0){\includegraphics{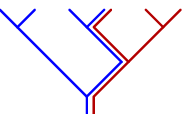}}
   \put(0,31){\sn{2}} \put(10,31){\sn{5}} \put(20,31){\sn{1}} 
   \put(32,31){\sn{4}}\put(42,31){\sn{3}}
   \put(31,49){\vector(0,-1){15}}
  \end{picture}
   \ \raisebox{12pt}{$\xrightarrow{\ \curlyvee\ }$}\quad
 \raisebox{7.5pt}{$\displaystyle  
   \raisebox{8pt}{$\biggl(\;$}
   \begin{picture}(32,27)
    \put(0,0){\includegraphics{figures/231.d.eps}}
    \put(0,21){\sn{2}} \put(10,21){\sn{5}} \put(20,21){\sn{1}} 
   \end{picture}
     ,\ 
   \begin{picture}(22,23)
    \put(0,0){\includegraphics{figures/21.d.eps}}
    \put(0,16){\sn{4}} \put(10,16){\sn{3}}
   \end{picture}
   \raisebox{8pt}{$\;\biggr)$}
   \quad \raisebox{5.5pt}{\mbox{or}}\quad 
   \raisebox{8pt}{$\biggl(\;$}
   \begin{picture}(32,27)
    \put(0,0){\includegraphics{figures/231.d.eps}}
    \put(0,21){\sn{2}} \put(10,21){\sn{3}} \put(20,21){\sn{1}} 
   \end{picture}
     ,\ 
   \begin{picture}(22,23)
    \put(0,0){\includegraphics{figures/21.d.eps}}
    \put(0,16){\sn{2}} \put(10,16){\sn{1}}
   \end{picture}
   \raisebox{8pt}{$\;\biggr)$}.
    $}
\]
Write $w\psplit(w_0,w_1)$ to indicate that the ordered forest $(w_0,w_1)$ 
(or pair of ordered trees) is obtained by splitting $w$ along some leaf. 
(Context will determine how to interpret the result.) 
More generally, we may split an ordered tree $w$ along a
multiset of $m\geq 0$ of its leaves to obtain an ordered forest, or tuple of ordered trees, written 
$w\psplit(w_0,\dotsc,w_m)$.
For example,
 \begin{equation}\label{Eq:splitting}
  \raisebox{-18pt}{$\displaystyle
   \begin{picture}(72,59)
    \put(0,0){\includegraphics{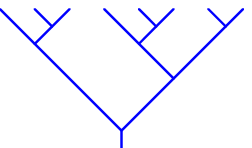}}
    \put( 0,41){\sn{3}} \put(10,41){\sn{2}} \put(21,41){\sn{7}} \put(30,41){\sn{5}} 
    \put(40,41){\sn{1}} \put(50,41){\sn{6}} \put(60,41){\sn{4}} 
    \put(19,59){\vector(0,-1){15}} \put(21,59){\vector(0,-1){15}}
    \put(50,59){\vector(0,-1){15}} \put(60,59){\vector(0,-1){15}}
   \end{picture}
    \ \raisebox{12pt}{\large$\leadsto$}\ 
   \begin{picture}(72,59)
    \put(0,0){\includegraphics{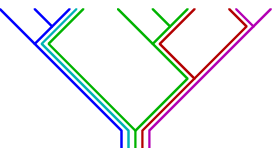}}
    \put( 0,41){\sn{3}} \put(10,41){\sn{2}} \put(24,41){\sn{7}} \put(34,41){\sn{5}} 
    \put(44,41){\sn{1}} \put(56,41){\sn{6}} \put(68,41){\sn{4}} 
    \put(21,59){\vector(0,-1){15}} \put(23,59){\vector(0,-1){15}}
    \put(55,59){\vector(0,-1){15}} \put(67,59){\vector(0,-1){15}}
   \end{picture}
    \quad\raisebox{12pt}{$\xrightarrow{\ \curlyvee\ }$}\quad
  \raisebox{5pt}{$\displaystyle 
   \raisebox{8pt}{$\biggl(\;$}
   \begin{picture}(22,23)
    \put(0,0){\includegraphics{figures/21.d.eps}}
    \put(0,16){\sn{3}} \put(10,16){\sn{2}}
   \end{picture}
     , \ 
     \includegraphics{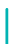}
    \ , 
   \begin{picture}(32,28)
    \put(0,0){\includegraphics{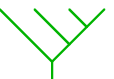}}
    \put(0,21){\sn{7}} \put(10,21){\sn{5}}\put(20,21){\sn{1}}
   \end{picture}
     , \ 
   \begin{picture}(12,18)
    \put(0,0){\includegraphics{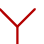}}
    \put(0,11){\sn{6}}
   \end{picture}
    \ , \ 
   \begin{picture}(12,18)
    \put(0,0){\includegraphics{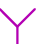}}
    \put(0,11){\sn{4}}
   \end{picture}
   \raisebox{8pt}{$\;\biggr)$}\ .$}$}
 \end{equation}

Given $v\in\s_m$ and an ordered forest $(w_0,\dotsc,w_m)$,
let \Blue{$(w_0,\dotsc,w_m)/v$} be the ordered tree
obtained by grafting the root of $w_i$ to the $i$th leaf of $v$, where the
nodes of $v$ are greater than all nodes of $w$, but the relative orders
within the $w_i$ and $v$ are maintained.
When $v$ is the ordered tree corresponding to $1432$ and
$w\psplit(w_0,\dotsc,w_m)$ is the splitting~\eqref{Eq:splitting}, this grafting is
\[
  \begin{picture}(125,88)
    \put(0,0){\includegraphics{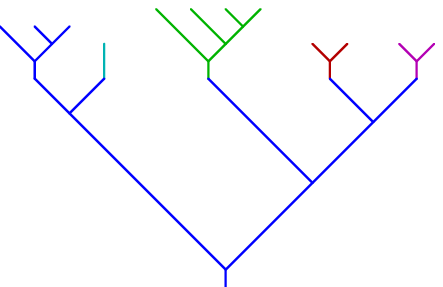}}
    \put(0,76){\sn{3}} \put(10,76){\sn{2}} \put(15,61){\sn{8}} \put(36,61){\sn{11}}
    \put(45,81){\sn{7}} \put(55,81){\sn{5}} \put(65,81){\sn{1}} \put(71,61){\sn{10}}
    \put(90,71){\sn{6}} \put(103,61){\sn{9}} \put(115,71){\sn{4}}
  \end{picture}
  \quad\raisebox{30pt}{$=$}\quad
  \begin{picture}(120,73)
    \put(0,0){\includegraphics{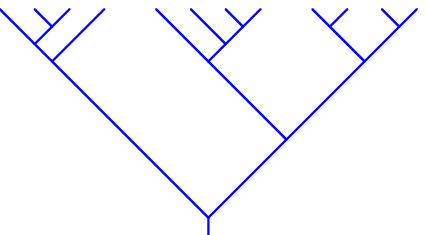}}
    \put(0,66){\sn{3}} \put(10,66){\sn{2}} \put(20,66){\sn{8}} \put(30,66){\sn{11}}
    \put(45,66){\sn{7}} \put(55,66){\sn{5}} \put(65,66){\sn{1}} \put(75,66){\sn{10}}
    \put(90,66){\sn{6}} \put(101,66){\sn{9}} \put(110,66){\sn{4}}
  \end{picture}
\]

The notions of splitting and grafting also make sense for the unordered trees $\y_n$ and we use the same notation, $\placeholder \psplit \placeholder$ and $\placeholder/\placeholder$. 
(Simply delete the labels in the constructions above.) 
These operations of splitting and grafting are compatible with the map 
$\tau\colon \s_\bb\to\y_\bb$:
if $w\psplit(w_0,\dotsc,w_m)$ then $\tau(w)\psplit(\tau(w_0),\dotsc,\tau(w_m))$ and all 
splittings in $\y_\bb$ are induced in this way from splittings in $\s_\bb$.
The same is true for grafting,
$\tau((w_0,\dotsc,w_m)/v)=(\tau(w_0),\dotsc,\tau(w_m))/\tau(v)$. 

%
%
\subsection{$\ssym$ and $\ysym$}\label{sec: ssym and ysym} 

For basics on Hopf algebras, see~\cite{Mont:1993}.
Let $\Blue{\ssym}:=\bigoplus_{n\geq0} {\ssym}_{n}$ be the graded $\Q$--vector space whose
$n^{\rm th}$ graded piece has  basis $\{F_w \mid w\in \s_n\}$. 
Malvenuto and Reutenauer~\cite{MalReu:1995} defined a Hopf algebra structure on
$\ssym$.
For $w\in\s_\bb$, define the coproduct
\[
   \Delta F_w\ :=\ \sum_{w\psplit (w_0,w_1)} F_{w_0}\otimes F_{w_1}\,,
\]
where $(w_0,w_1)$ is a pair of ordered trees. If $v\in\s_m$, define the product
\[
  F_w \cdot F_v\ :=\ \sum_{w\psplit(w_0,\dotsc,w_m)} F_{(w_0,\dotsc,w_m)/v}\,.
\]
The counit is the projection $\varepsilon\colon\ssym \to \ssym_0$ onto the $0$th graded piece,  which is spanned by the unit, $1=F_\emptyset$, for this multiplication.

\begin{prop}[\cite{MalReu:1995}]
  With these definitions of coproduct, product, counit, and unit, $\ssym$ is a 
  graded, connected cofree Hopf algebra that is neither commutative nor cocommutative.
\end{prop}

Let $\Blue{\ysym}:=\bigoplus_{n\geq0} {\ysym}_{n}$ be the graded $\Q$--vector space whose
$n^{\rm th}$ graded piece has  basis $\{F_t \mid t\in \y_n\}$. 
Loday and Ronco~\cite{LodRon:1998} defined a Hopf algebra structure on
$\ysym$.
For $t\in\y_\bb$, define the coproduct
\[
   \Delta F_t\ :=\ \sum_{t\psplit (t_0,t_1)} F_{t_0}\otimes F_{t_1}\,,
\]
and if $s\in\y_m$, define the product
\[
  F_t \cdot F_s\ :=\ \sum_{t\psplit(t_0,\dotsc,t_m)} F_{(t_0,\dotsc,t_m)/s}\,.
\]
The counit is the projection $\varepsilon\colon\ysym \to \ysym_0$ onto the $0$th graded piece, which is
spanned by the unit, $1=F_{\includegraphics{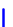}}$, for this multiplication.
The map $\tau$ extends to a linear map $\btau\colon\ssym\to\ysym$, defined by
$\btau(F_w) = F_{\tau(w)}$.

\begin{prop}[\cite{LodRon:1998}]
  With these definitions of coproduct, product, counit, and unit, $\ysym$ is a 
  graded, connected cofree Hopf algebra that is neither commutative nor cocommutative
  and the map $\btau$ a morphism of Hopf algebras.
\end{prop}

Some structures of the Hopf algebras $\ssym$ and $\ysym$, particularly their primitive elements
and coradical filtrations are better understood with respect to a second basis.
The M\"obius function $\mu$ (or $\mu_P$) of a poset $P$ is defined for pairs $(x,y)$ of elements
of $P$ with $\mu(x,y)=0$ if $x\not< y$, $\mu(x,x)=1$, and, if $x<y$, then 
 \begin{equation}\label{Eq:Moebius_def}
   \mu(x,y)\ =\ -\sum_{x\leq z<y}\mu(x,z)
   \qquad\mbox{so that}\qquad
   0\ =\ \sum_{x\leq z\leq y}\mu(x,z)\ .
 \end{equation}

For $w\in\s_\bb$ and $t\in\y_\bb$, set
\begin{equation}\label{eq: ssym ysym M-basis}
  \Blue{M_w}\ :=\ \sum_{w\leq v} \mu(w,v) F_v
  \qquad\mbox{and}\qquad
  \Blue{M_t}\ :=\ \sum_{t\leq s} \mu(t,s) F_s\ ,
\end{equation}
where the first sum is over $v\in\s_\bb$, the second sum over $s\in\y_\bb$, and 
$\mu(\cdot, \cdot)$ is the M\"obius function in the weak and Tamari orders.

\begin{prop}[\cite{AguSot:2005,AguSot:2006}]\label{prop:M-basis}
If $w\in\s_\bb$, then
\begin{equation}\label{eq:tau_M-basis}
  \btau(M_w)\ =\ 
    \begin{cases} M_{\tau(w)}, & \mbox{if }\,w=\Max[\tau(w)]\,,\\
                        0, & \mbox{otherwise}\,
    \end{cases}
\end{equation}
 and  
 \begin{equation}\label{Eq:SSym_M_coprod}
   \Delta M_w\ =\  \sum_{w=u\backslash v} M_u\otimes M_v\,.
 \end{equation} 
 If $t\in\y_\bb$, then
 \begin{equation} \label{Eq:YSym_M_coprod}
   \Delta M_t\ =\  \sum_{t=r\backslash s} M_r\otimes M_s\,.
 \end{equation} 
\end{prop}

This implies that the set $\{M_w\mid w\in\s_\bb\mbox{ is indecomposable}\}$ is a basis 
for the primitive elements of $\ssym$ (and the same for $\ysym$), 
thereby explicitly realizing the cofree-ness of $\ssym$ and $\ysym$.

\subsection{M\"obius functions and interval retracts}\label{sec: poset lemma}

A pair $f\colon P\to Q$ and $g\colon Q\to P$ of poset maps is a \demph{Galois connection} if
$f$ is left adjoint to $g$ in that 
\[
   \forall\;p\in P\mbox{ and }q\in Q\,, \quad
     f(p)\ \leq_Q\ q\ \Longleftrightarrow\ p\ \leq_P\ g(q)\,.
\]
When this occurs, Rota~\cite[Theorem 1]{Rot:1964} related the M\"obius functions of $P$ and
$Q$:
\[
   \forall\;p\in P\mbox{ and }q\in Q\,, \quad
   \sum_{f(y)=q} \mu_P(p,y)\ \ =\ \, 
   \sum_{g(x)=q} \mu_Q(x,q)\,.
\]
Rota's formula was used in \cite{AguSot:2006} to establish the coproduct formulas~\eqref{eq:tau_M-basis}
and~\eqref{Eq:YSym_M_coprod}, as the maps $\tau\colon\s_\bb\to\y_\bb$ and
$\Max\colon\y_\bb\to\s_\bb$ form a Galois connection~\cite[Section 9]{BjoWac:1997}.

We do not have a Galois connection between $\s_\bb$ and $\m_\bb$, and so cannot use Rota's
formula. 
Nevertheless, there is a useful relation between the M\"obius functions of 
$\s_\bb$ and $\m_\bb$ that we establish here in a general form. 
A surjective poset map $f\colon P\to Q$ from a finite lattice $P$ is an 
\demph{interval retract} if the fibers of $f$ are intervals and if $f$ admits an
order-preserving section $g\colon Q\to P$ with $f\circ g=\mbox{id}$.

\begin{theorem}\label{thm: subGalois} 
Let the poset map $f\colon P\to Q$ is an interval retract, then the
M\"obius functions $\mu_P$ and $\mu_Q$ of $P$ and $Q$ are related by the formula 
\begin{equation}\label{eq: subGalois}
    \mu_Q(x,y) \ =\,  
    \sum_{\substack{f(a)=x\\f(b)=y}} \mu_P(a,b) \qquad (\forall x,y\in Q).
\end{equation}
\end{theorem}

In Section~\ref{sec: m}, we define an interval retract $\beta\colon \s_n\to\m_n$.

We evaluate each side of \eqref{eq: subGalois} using Hall's formula, which expresses the
M\"obius function in terms of \emph{chains}.  
A linearly ordered subset $C\colon x_0<\dotsb<x_r$ of a poset is a \demph{chain} of
\demph{length $\ell(C)=r$} from $x_0$ to $x_r$. 
Given a poset $P$, let $\Chains(P)$ be the set of all chains in $P$.
A poset $P$ is an \demph{interval} if it has a unique maximum element and a unique minimum
element.
If $P=[x,y]$ is an interval, let $\Chains'(P)$ denote the chains in $P$ beginning in $x$ and
ending in $y$.  
Hall's formula states that 
\[
	\mu(x,y) \ \,=\  \sum_{C \in \chains'[x,y]} \!(-1)^{\ell(C)} \,.
\]
Our proof rests on the following two lemmas.

\begin{lemma}\label{thm: interval} If $P$ is an interval, then 
$\displaystyle
	\sum_{C \in \chains(P)} (-1)^{\ell(C)} = 1.
$ 
\end{lemma}

\begin{proof}
 Suppose that $P=[x,y]$ and append new minimum and maximum elements to P to get 
 $\hat P := P \cup \{\hat0, \hat1\}$.  
 Then the definition of M\"obius function~\eqref{Eq:Moebius_def} gives
\[
   \mu(\hat0,\hat1)\ =\ -\sum_{\hat0\leq z \leq y}\mu(\hat0,z)\,,
\]
 which is zero by~\eqref{Eq:Moebius_def}.
 By Hall's formula,
\[
  0\ =\ \mu(\hat0,\hat1)\ =\  \sum_{C \in \chains'[\hat0,\hat1]} (-1)^{\ell(C)}\ =\ 
  -1 + \sum_{C\in\chains(P)} (-1)^{\ell(C)+2}\,,
\]
 where the term $-1$ comes from the chain $\hat0<\hat1$.
 This proves the lemma.
\end{proof}

Call a partition $P=K_0\sqcup \dotsb \sqcup K_r$ of $P$ into subposets $K_i$  \demph{monotone} if 
$x<y$ with $x\in K_i$ and $y\in K_j$ implies that $i\leq j$. 
Given $\emptyset \subsetneq I \subseteq [0,r]$, write $\Chains_I(P)$ for the
subset of chains $C$ in $\Chains(P)$ such that $C \cap K_i \neq \emptyset$
if and only if $i\in I$. 

\begin{lemma}\label{thm: chain preimages} 
Let $P=K_0\sqcup \dotsb \sqcup K_r$ be a monotonic partition of a poset $P$. 
If $\bigcup_{i\in I} K_i$ is an interval for all $I\subseteq[0,r]$, then 
\begin{equation}\label{eq: monotonic chains}
	\sum_{C \in \chains_{[0,r]}(P)} (-1)^{\ell(C)} \ =\ (-1)^r \,.
\end{equation}
\end{lemma}

\begin{proof}
We argue by induction on $r$. 
Lemma~\ref{thm: interval} is the case $r=0$ (wherein $K_0=P$), so we consider the case $r\geq1$.

Form the poset $\hat P = P \cup \{\hat0,\hat1\}$ as in the proof of Lemma \ref{thm: interval}.
Since $P$ is an interval, we have
$
\sum_{C \in \chains'[\hat0,\hat1]} (-1)^{\ell(C)} = 0.
$ 
As $\Chains'[\hat0,\hat1]=\bigsqcup_I\Chains_I(P)$ we have,
\[
    0 \ =\ -1 \ +\  
      \sum_{\emptyset\subsetneq I \subsetneq [0,r]} 
            \biggl(\,\sum_{C \in \chains_I(P)} (-1)^{\ell(C)} \biggr)
      \ +\  \sum_{C \in \chains_{[0,r]}(P)} (-1)^{\ell(C)} \,,
\] 
where the term $-1$ counts the chain $\hat0<\hat1$. 
Applying induction, we have
\[
0 = \sum_{k=0}^{r} \binom{r+1}{k} (-1)^{k-1} + \sum_{C\in \chains_{[0,r]}(P)}(-1)^{\ell(C)} \,.
\]
Comparing this to the binomial expansion of $(1-1)^{r+1}$ completes the proof.
\end{proof}

\begin{proof}[Proof of Theorem \ref{thm: subGalois}]
Fix $x<y$ in $Q$. 
We use Hall's formula to rewrite the right-hand side of \eqref{eq: subGalois} as
 \begin{equation}\label{Eq:first_expansion}
	\sum_{\substack{f(x)=a\\f(y)=b}}\; \sum_{C\in
          \chains'[a,b]}(-1)^{\ell(C)}\,. 
 \end{equation}
Fix a chain $D\colon q_0<\dotsb<q_r$ in $\Chains'[x,y]$ and 
let \Blue{$P|_D$} be the subposet of
$P$ consisting of elements that occur in some chain of $P$ that maps to $D$ under $f$.
This is nonempty as $f$ has section.
Furthermore, the sets $\Blue{K_i}:=f^{-1}(q_i)\cap P|_D$, for $i=0,\dotsc,r$, form a
monotonic partition of $P|_D$.
We claim that $\bigcup_{i\in I} K_i$ is an interval for all $I\subseteq[0,r]$.
If so, let us first rewrite~\eqref{Eq:first_expansion} as a sum over chains $D$ in $Q$,
\[
	\sum_{D\in\chains'[x,y]} \ 
	\sum_{C\in\chains_{[0,\ell(D)]}(P|_D)} (-1)^{\ell(C)} \,.
\]
By Lemma~\ref{thm: chain preimages}, the inner sum becomes
$\sum_{D}(-1)^{\ell(D)}$, which completes the proof.

To prove the claim, suppose that $I = \{i_0<\dotsb<i_s\}$.
Each set $K_i$ ($i\in I$) is an interval, as it is the intersection of two
intervals in the lattice $P$. 
Thus $K_{i_0}$ and $K_{i_s}$ are intervals with minimum and
maximum elements $m$ and $M$, respectively.
Any chain in $\bigcup_{i\in I} K_i$ can be extended to a chain beginning with $m$ and ending
at $M$, so $\bigcup_{i\in I} K_i$ is an interval.
\end{proof}

\section{The Multiplihedra $\m_\bb$}\label{sec: m}

The map $\tau\colon\s_\bb\to\y_\bb$ forgets the linear ordering
of the node poset of an ordered tree, and it induces a morphism of Hopf algebras
$\btau\colon\ssym\to\ysym$.
In fact, one may take the (ahistorical) view that the Hopf structure on $\ysym$ is induced
from that on $\ssym$ via the map $\tau$.
Forgetting some, but not all, of the structure on a tree in $\s_\bb$ factorizes the map $\tau$.
Here, we study  combinatorial consequences of one such factorization, and later treat
its algebraic consequences.

\subsection{Bi-leveled trees}\label{sec: bi-leveled trees}

A \demph{bi-leveled tree} $(t;\setT)$ is a planar binary tree $t\in \y_n$ together with 
an (upper) order ideal $\setT$ of its node poset, where $\setT$ contains the leftmost node
of $t$ as a minimal element.
Thus $\setT$ contains all nodes along the path from the leftmost leaf to the root, and none above the leftmost node.
Numbering the gaps between the leaves of $t$ by $1,\dotsc,n$ from left to
right, $\setT$ becomes a subset of $\{1,\dotsc,n\}$.

Saneblidze and Umble~\cite{SanUmb:2004} introduced bi-leveled trees 
to describe a cellular projection from the permutahedra to Stasheff's 
multiplihedra $\m_\bb$, with 
the bi-leveled trees on $n$ nodes indexing the vertices $\m_n$.
Stasheff used a different type of tree for the vertices of $\m_\bb$.
These alternative trees lead to a different Hopf structure which we explore in 
a forthcoming paper~\cite{FLS_CCC}.
We remark that $\m_0 = \{\raisebox{-1pt}{\includegraphics{figures/0.eps}}\}$. 

\begin{figure}[htb]
\[
  \begin{picture}(192,172)(0,2)
   \put(2,5){\includegraphics[height=170pt]{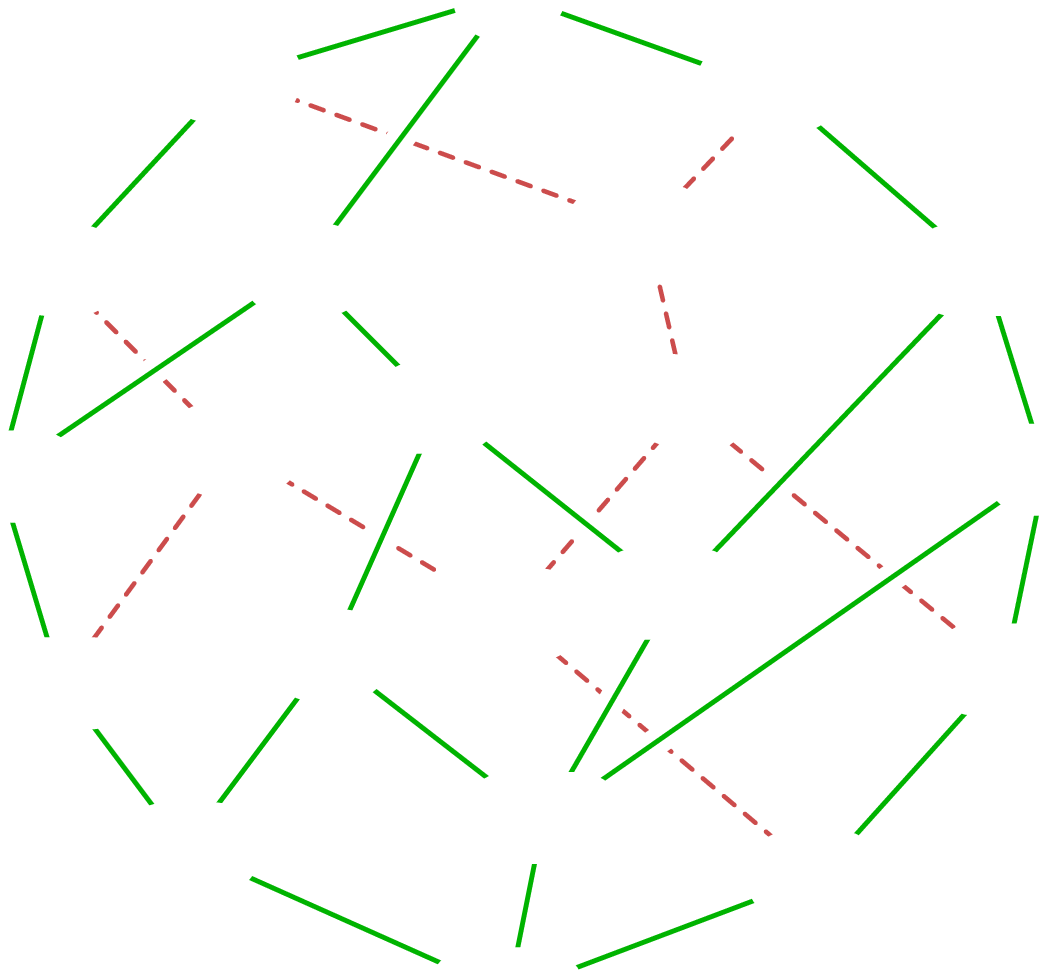}}

   \put(84,162){\includegraphics[height=12pt]{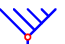}} 

   \put(39.5,149){\includegraphics[height=12pt]{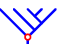}} 
   \put(127,146.5){\includegraphics[height=12pt]{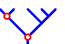}} 

   \put(11,118){\includegraphics[height=12pt]{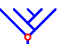}}
   \put(51,118){\includegraphics[height=12pt]{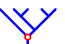}}
   \put(106.5,124.5){\includegraphics[height=12pt]{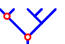}}
   \put(160.7,118){\includegraphics[height=12pt]{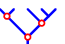}}

   \put( 2,84.5){\includegraphics[height=12pt]{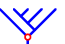}}
   \put(39.7,89){\includegraphics[height=12pt]{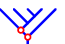}}
   \put(72.7,95.7){\includegraphics[height=12pt]{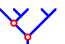}}
   \put(112.8,95.5){\includegraphics[height=12pt]{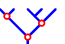}}
   \put(170.3,85.7){\includegraphics[height=12pt]{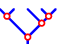}}

   \put(12, 51.5){\includegraphics[height=12pt]{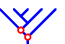}}
   \put(55.3, 55.7){\includegraphics[height=12pt]{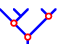}}
   \put(84, 62.5){\includegraphics[height=12pt]{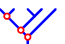}}
   \put(110.5, 66){\includegraphics[height=12pt]{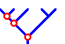}}
   \put(163.8, 53.5){\includegraphics[height=12pt]{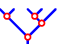}}

   \put(32, 25){\includegraphics[height=12pt]{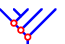}}
   \put(89.5, 29.5){\includegraphics[height=12pt]{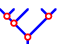}}
   \put(134, 20.3){\includegraphics[height=12pt]{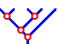}}

   \put(84, 1.5){\includegraphics[height=12pt]{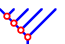}}

   \end{picture}
\]

\caption{The $1$-skeleton of the multiplihedron $\m_4$.}
\label{fig: M4}
\end{figure}

The partial order on $\m_n$ is defined by 
$(s;\setS)\leq(t;\setT)$ if $s\leq t$ in $\y_n$ and $\setS\supseteq \setT$. 
The Hasse diagrams of the posets $\m_n$ are $1$-skeleta for the multiplihedra.
We represent a bi-leveled tree by drawing the underlying tree $t$ and circling the nodes in
$\setT$. 
The Hasse diagram of $\m_4$ appears in Figure \ref{fig: M4}. 

\subsection{Poset maps}\label{sec: maps}

Forgetting the order ideal in a bi-leveled tree, $(t;\setT)\mapsto t$, is a poset map 
$\phi\colon\m_\bb\to\y_\bb$.
We define a map $\beta\colon\s_\bb\to\m_\bb$ so that
\[
   \s_\bb\ \xrightarrow{\ \beta\ }\ \m_\bb
         \ \xrightarrow{\ \phi\ }\ \y_\bb
\]
factors the map $\tau\colon\s_\bb\to\y_\bb$, and we define a right inverse (section) $\iota$ of 
$\beta$. 

Let $w\in\s_\bb$ be an ordered tree. 
Define the set
 \begin{equation}\label{Eq:T(w)}
   \setT(w)\ :=\ \{i\mid w(i)\geq w(1)\}\,.
 \end{equation}
Observe that $(\tau(w);\setT(w))$ is a bi-leveled tree.
Indeed, as $w$ is a linear extension of $\tau(w)$, $\setT(w)$ is an upper order ideal which by
definition~\eqref{Eq:T(w)} contains the leftmost node as a minimal element.
Since covers in the weak order can only decrease the subset $\setT(w)$ and $\tau$ is also a
poset map, we see that $\beta$ is a poset map.

\begin{theorem}
  The maps $\beta\colon\s_\bb\to\m_\bb$ and $\phi\colon\m_\bb\to\y_\bb$ are surjective poset
  maps with $\tau=\phi\circ\beta$.
\end{theorem}

\begin{figure}[!ht]
\[
  \raisebox{40pt}{$
     \beta^{-1}\biggl(\;\raisebox{-15pt}{\includegraphics{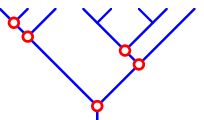}}\;\biggr)
     \quad =$}\quad
  \begin{picture}(145,93)(-25,0)
   \put(-5,3){\includegraphics{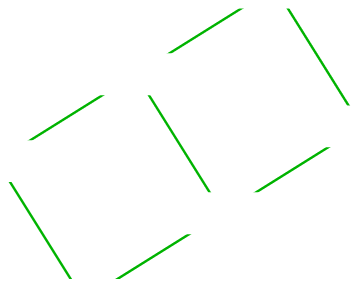}}
   \put(0,0){$3471526$}
   \put(40,25){$3571426$}
   \put(80,50){$3671425$}
   \put(-25,40){$3472516$}
   \put(15,65){$3572416$}
   \put(55,90){$3672415$}
  \end{picture}
\]
\caption{The preimages of $\beta$ are intervals.}
\label{fig: preimages are intervals}
\end{figure}

\begin{figure}[!ht]
\[
  \begin{picture}(400,144)(0,-12)
    \put( 42,80){\begin{picture}(80,40)
       \put(0,0){\includegraphics{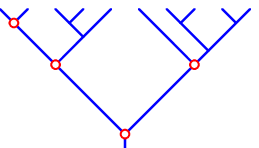}}
       \put(-2,32){\Red{{\tiny 1}}}  \put(9,20){\Red{{\tiny 2}}}
       \put(13,32){{\tiny 3}}  \put(26,28){{\tiny 4}}
       \put(34,10){\Red{{\tiny 5}}}  \put(48,21){\Red{{\tiny 6}}}
       \put(56,34){{\tiny 7}} \put(62,24){{\tiny 8}} 
       \put(70,32){{\tiny 9}} 
     \end{picture}}
    \put(138,100){\vector(1,0){20}}  \put(140,100){\vector(-1,0){20}}
    \put(164,80){\includegraphics{figures/m674598231.s.eps}}
    \put(260,100){\vector(1,0){20}}  \put(262,100){\vector(-1,0){20}}
    \put(286,80){\includegraphics{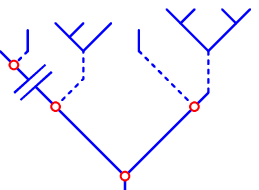}}
    \put(0,20){$\left(\,\,\raisebox{-17pt}{\includegraphics{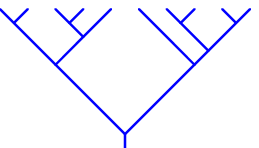}}
                \!;\,\{1,2,5,6\}\right)$}
    \put( 78,54){\vector(0,1){20}}  \put( 78,70){\vector(0,-1){20}}
    \put(322,54){\vector(0,1){20}}  \put(322,66){\vector(0,-1){20}}
    \put(243,20){$\left(\,\,\raisebox{-10pt}{\includegraphics{figures/231.d.eps}}
            \!,\, \left(\,\raisebox{-4pt}{\includegraphics{figures/0.d.eps}}
              ,\raisebox{-6pt}{\includegraphics{figures/12.d.eps}}
              ,\raisebox{-4pt}{\includegraphics{figures/0.d.eps}}
              ,\raisebox{-8pt}{\includegraphics{figures/231.d.eps}}\right)\right)$}
    \put( 65,-12){$(t;\setT)$}
    \put(309,-12){$(t_0,\bm t)$}
  \end{picture}
\] 
\caption{Two representations of bi-leveled trees.}
\label{fig:correspondence}
\end{figure}

The fibers of the map $\beta$ are intervals (indeed, products of intervals); see
Figure~\ref{fig: preimages are intervals}. 
We prove this using an equivalent representation of a bi-leveled tree and a description of
the map $\beta$ in that representation.
If we prune a bi-leveled tree $b=(t;\setT)$ above
the nodes in $\setT$ (but not on the leftmost branch) we obtain a tree $t_0'$ (the order
ideal) on $r$ nodes and a planar forest $\bm t= (t_1,\ldots,t_r)$ of $r$ trees. 
If we prune $t_0'$ just below its leftmost node, we obtain the tree
$\includegraphics{figures/1.eps}$ (from the pruning) and a tree $t_0$, 
and $t_0'$ is obtained by grafting 
$\includegraphics{figures/1.eps}$ onto the leftmost leaf of  $t_0$.
We may recover $b$ from this tree $t_0$ on $r{-}1$ nodes and the planar
forest $\bm t= (t_1,\ldots,t_r)$, and so we also write
$b=(t_0,\bm t)$.
We illustrate this correspondence in Figure~\ref{fig:correspondence}.
%

We describe the map $\beta$ in terms of this second representation of bi-leveled trees.
Given a permutation $w$ with $\beta(w)=(t;\setT)$ and $|\setT|=r$, let
$u_1u_2\dotsc u_r$ be the restriction of $w$ to the set $\setT$.
We may write the values of $w$ as $w=u_1v^1u_2\dotsb u_rv^r$, where $v^i$ is the
(possibly empty) subword of $w$ between the numbers $u_i$ and $u_{i+1}$ and $v^r$ is the
word after $u_r$.
Call this the \demph{bi-leveled factorization} of $w$.
For example, 
\newcommand{\fstack}[2]{%
	\begin{array}{@{}c@{}} {\textstyle #1}\\[-.25ex] {\scriptstyle #2} \end{array}
}
\[
  \Blue{6}54\Blue{789}132\ \longmapsto\ (
	\raisebox{-1.25ex}{$\fstack{\Blue{6}\scriptscriptstyle\,}{\red u_1}$},\,
	\raisebox{-1.25ex}{$\fstack{54}{\red v^1}$}\,,\,
	\raisebox{-1.25ex}{$\fstack{\Blue{7}}{\red u_2}$},\,
	\raisebox{-1.25ex}{$\fstack{\emptyset\,}{\red v^2}$},\,
	\raisebox{-1.25ex}{$\fstack{\Blue{8}}{\red u_3}$},\,
	\raisebox{-1.25ex}{$\fstack{\emptyset\,}{\red v^3}$},\,
	\raisebox{-1.25ex}{$\fstack{\Blue{9}}{\red u_4}$},\,
	\raisebox{-1.25ex}{$\fstack{132}{\,\red v^4}$}{\scriptscriptstyle\,}
	)\,.
\]
Note that 
$\beta(w)=(\tau(\overline{u_2\dotsc u_r}), (\tau(\overline{v^1}),\dotsc,\tau(\overline{v^r})))$.

\begin{theorem}\label{thm:beta_preimages}
  For any $b \in \m_n$ the fiber $\beta^{-1}(b) \subseteq \s_n$ is a
  product of intervals. 
\end{theorem}

\begin{proof}
 Let $b=(t_0,(t_1,\dotsc,t_r))=(t;\setT)\in\m_n$ be a bi-leveled tree.
 A permutation $w\in\beta^{-1}(b)\in\s_n$ has a bi-leveled factorization
 $w=u_1v^1u_2\dotsc u_rv^r$ with 
 \begin{equation}\label{eq:bi-leveled_factorization}
  \begin{array}{ll}
    (\mbox{{\it i}})&w|_\setT=u_1u_2\dotsc u_r,\   u_1=n{+}1{-}r, \ 
        \tau(\overline{u_2\dotsc u_r})=t_0,\ \mbox{ and }\\
    (\mbox{{\it ii}})&\tau(\overline{v^i})=t_i, \mbox{ for } i=1,\dotsc,r\,.\rule{0pt}{14pt}
   \end{array}
 \end{equation}
 Since $u_1<u_2,\dotsc,u_r$ are the values of $w$ in the positions of $\setT$,
 and $u_1=n{+}1{-}r$ exceeds all the letters in $v^1,\dotsc,v^r$, which are the values of $w$ in
 the positions in the complement of $\setT$, these two parts of the bi-leveled factorization
 may be chosen independently to satisfy~\eqref{eq:bi-leveled_factorization}, which shows that
 $\beta^{-1}(b)$ is a product.

 To see that the factors are intervals, and thus $\beta^{-1}(b)$ is an interval, we examine
 the conditions $(\mbox{\it i})$ and $(\mbox{\it ii})$ separately.
 Those $u_1\dotsc u_r=w|_\setT$ for $w$ in the fiber $\beta^{-1}(b)$ are exactly the
 set of $n{+}1{-}r,u_2,\dotsc,u_r$ with $\{u_2,\dotsc,u_r\}=\{n{+}2{-}r,\dotsc,n\}$ and
 $\tau(\overline{u_2\dotsc u_r})=t_0$.
 This is a poset under the restriction of the weak order, and it  
 is in natural bijection with the interval $\tau^{-1}(t_0)\subset\s_{r-1}$.
 Its minimal element is $\Blue{\Min_0(b)}=u_1u_2\dotsc u_r$, where $u_2\dotsc u_r$ is the
 unique $231$-avoiding word on  $\{n{+}1{-}r,\dotsc,n\}$ satisfying 
 $(\mbox{\it i})$, and its maximal element is
 $\Blue{\Max_0(b)}=u_1u_2\dotsc u_r$, where now $u_2\dotsc u_r$ is the unique $132$-avoiding
 word on $\{n{+}1{-}r,\dotsc,n\}$ satisfying $(\mbox{\it i})$.  

 Now consider sequences of words $v^1,\dotsc,v^r$ on distinct letters $\{1,\dotsc,n{-}r\}$
 satisfying $(\mbox{\it ii})$.
 This is also a poset under the restriction of the weak order.
 It has a minimal element, which is the unique such sequence \Blue{$\underline{\Min}(b)$} 
 satisfying $(\mbox{\it ii})$ where the letters of $v^i$ preceed those of $v^j$ whenever $i<j$,
 and where each $v^i$ is $231$-avoiding.
 Its maximal element is the unique sequence \Blue{$\underline{\Max}(b)$} 
 satisfying $(\mbox{\it ii})$ where the letters of $v^i$ are greater than those of $v^j$ when
 $i<j$ and $v^i$ is $132$-avoiding. 
\end{proof}

\begin{figure}[htb]
\[
 \begin{picture}(150,85)
   \put(70,34){\vector(1,0){30}}\put(80,38){$\beta$}

   \thicklines
   \put(0,70){$2143$}
   \put(110,60){\includegraphics{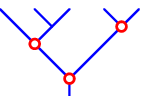}}

   \put(129,30){\Color{1 0 1 .3}{\line(0,1){25}}}
   \put(13,12){\Color{1 0 1 .3}{\line(0,1){52}}}
   \put(23,12){\Color{1 0 1 .3}{\line(2,1){15}}}
   \put(0,0){$1243$}  \put(40,20){$1342$}
   \put(110, 0){\includegraphics{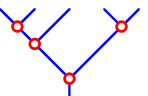}}
   \end{picture}
  \qquad\qquad\quad 
  \begin{picture}(150,85)
   \put(70,34){\vector(1,0){30}}\put(80,38){$\beta$}

   \thicklines
   \put(0,70){$3241$}  \put(40,50){$3142$}
   \put(110,60){\includegraphics{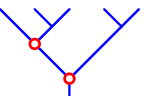}}

   \put(129,30){\Color{1 0 1 .3}{\line(0,1){25}}}
   \put(13,12){\Color{1 0 1 .3}{\line(0,1){52}}}
   \put(0,0){$2341$}  
   \put(23,65){\Color{1 0 1 .3}{\line(2,-1){15}}}
   \put(110, 0){\includegraphics{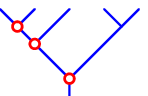}}
   \end{picture}
\]
\caption{Fibers of $\beta$.}
\label{Fig:not_order}
\end{figure}

The fibers of $\beta$ are intervals so that consistently choosing the minimum or maximum
in a fiber gives two set-theoretic sections.
These are not order-preserving as may be seen from Figure~\ref{Fig:not_order}.
We have $\raisebox{-2pt}{\includegraphics{figures/m2143.s.eps}}<
   \raisebox{-2pt}{\includegraphics{figures/m1243.s.eps}}$ but the maxima
   in their fibers under $\beta$, $1342$ and $2143$, are incomparable.
Similarly, $\raisebox{-2pt}{\includegraphics{figures/m3241.s.eps}}<
   \raisebox{-2pt}{\includegraphics{figures/m2341.s.eps}}$ but the minima
   in their fibers under $\beta$, $2341$ and $3142$, are incomparable.
This shows that the map $\beta\colon\s_\bb\to\m_\bb$ is not a lattice congruence (unlike the map $\tau \colon \s_\bb \to \y_\bb$~\cite{Re04}). 

 In the notation of the proof, given a bi-leveled tree $b=(t_0,(t_1,\dotsc,t_r))$, 
 let $\Blue{\iota(b)}$ be the permutation $w\in\beta^{-1}(b)$ with bi-leveled factorization 
 $w=u_1v^1u_2\dotsc u_rv^r$ where $u_1u_2\dotsc u_r=\Min_0(b)$  and 
 $(v^1,\dotsc,v^r)=\underline{\Max}(b)$.
 This defines a map $\iota\colon\m_n\to\s_n$ that is a section of the map $\beta$.
For example,
\[
   \iota
   \Biggl(\ \raisebox{-24pt}{\includegraphics{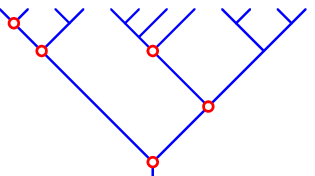}}\ \Biggr)\
  \ =\ 
  \raisebox{-32.5pt}{\begin{picture}(120,73)
    \put(0,0){\includegraphics{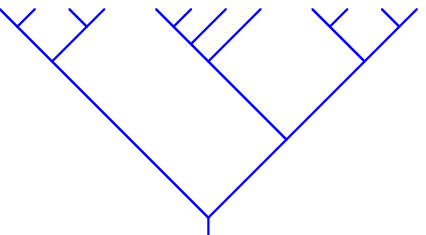}}
    \put(0,66){\sn{\Blue{7}}} \put(10,66){\sn{\Blue{8}}} \put(20,66){\sn{6}} 
    \put(30,66){\sn{\Blue{11}}}    \put(45,66){\sn{4}} \put(55,66){\sn{5}} 
    \put(65,66){\sn{\Blue{9}}} \put(75,66){\sn{\Blue{10}}}
    \put(90,66){\sn{2}} \put(101,66){\sn{3}} \put(110,66){\sn{1}}
  \end{picture}}
  =\ \Blue{78}6\,\Blue{11}\,45\Blue{9\,10\,}231\,.
\]

\begin{remark}\label{Rem:pattern}
 This map $\iota$ may be characterized in terms of pattern avoidance: the permutation $\iota(b)$ is the unique
 $w\in \beta^{-1}(b)$ avoiding the pinned patterns
 \[
	\bigl\{\underline{2}031, \,  \underline{0}231, \, \underline{3}021 \bigr\}\,,
 \]
where the underlined letter must be the first letter of a permutation. 
To see this, note that the first pattern forces the letters in
$v^i$ to be larger than those in $v^{i+1}$ for $1\leq i < r$, the
second pattern forces $u_2\dotsc u_r$ to be 231-avoiding, and the 
last pattern forces each $v^i$ to be 132-avoiding.
\end{remark}

\begin{theorem} \label{thm:iota}
  The map $\iota$ is injective, right-inverse to $\beta$, and order-preserving. 
  That is, $\beta\colon\s_n\to \m_n$ is an interval retract.
\end{theorem}


Since $\s_n$ is a lattice \cite{GuiRos:1960}, the fibers of $\beta$ are
intervals, and $\iota$ is a section of $\beta$. 
That is, we need only verify that $\iota$ is order-preserving. 
We begin by describing the covers in $\m_\bb$. 
Since $\beta$ is a surjective poset map, every cover in $\m_n$ is the image of 
some cover $w\lessdot w'$ in $\s_n$. 

\begin{lemma}\label{lem:covers}
If a cover $w\lessdot w'\in \s_n$ does not collapse under $\beta$, i.e., $\beta(w) \neq \beta(w')$, then it yields one of three types of covers $\beta(w)\lessdot \beta(w')$ in $\m_n$. 
\begin{enumerate}
 \item[$(i)$]
          In exactly one tree $t_i$ in $\beta(w)=(t_0,(t_1,\dotsc,t_r))$, a node is moved 
          from left to right across its parent to obtain $\beta(w')$. That is, $t_i \lessdot t'_i$.

 \item[$(ii)$]
          If $\beta(w)=(t;\setT)$, the leftmost node of $t$ is moved across its parent, 
          which has no other child in the order ideal $\setT$, and is deleted from $\setT$ to 
          obtain $\beta(w')$.

 \item[$(iii)$]
          If $\setT(w)=\{1=T_1<\dotsb<T_r\}$, then $\tau(w')=\tau(w)$ and $\setT(w')=\setT(w)\setminus\{T_j\}$ for some $j>2$.

\end{enumerate}
\end{lemma}

\begin{proof}
Put $w'=(k,k{+}1)w$, with $k,k{+}1$ appearing in order in $w$. 
Let $(t;\setT)$ and $(t_0,(t_1,\dotsc,t_r))$ be the two representations of $\beta(w)$. 
Write $\setT=\{T_1<\dotsb<T_r\}$ (with $T_1=1$) and $w|_\setT=u_1u_2\dotsc u_r$. 
If $w\lessdot w'$ and $\beta(w)\lessdot\beta(w')$, then $k$ appears 
within $w$ in one of three ways: 
	$(i)$ $u_1\neq k$, 
	$(ii)$ $u_1= k$ and $u_2=k{+}1$, or 
	$(iii)$ $u_1=k$ and $u_j=k{+}1$ for some $j>2$. 
These yield the corresponding descriptions in the statement of the lemma. 
(Note that in type $(i)$, $\setT(w')=\setT$, so if we set $\beta(w')=(t'_0,(t'_1,\dotsc,t'_r))$, 
then $t_i=t'_i$, except for one index $i$, where $t_i\lessdot t'_i$.) 
\end{proof}

\begin{figure}[htb]

\[
  \begin{picture}(340,120)
   \put(  0,80){\includegraphics{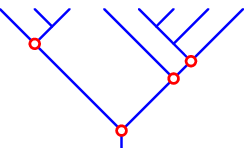}}
   \put( 90,80){\includegraphics{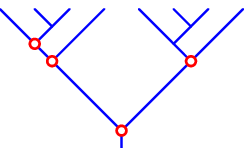}}
   \put(180,80){\includegraphics{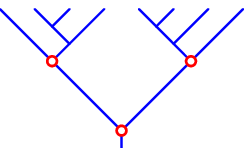}}
   \put(270,80){\includegraphics{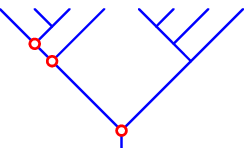}}

   \put( 55,50){$(i)$}
   \put(140,68){$(i)$}
   \put(180,68){$(ii)$}
   \put(265,50){$(iii)$}
   
   \put(35,45){\includegraphics{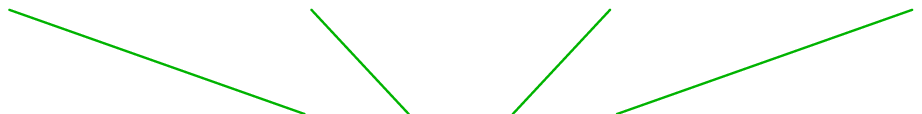}}

   \put(135,0){\includegraphics{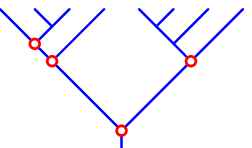}}
  \end{picture}
\]
\caption{Some covers in $\m_7$.}
\label{Fig:M-cover} 
\end{figure}

Figure~\ref{Fig:M-cover} illustrates these three types of covers, labeled
by their type.

For $\setT\subset\{1,\dotsc,n\}$ with $1\in\setT$, let
$\Blue{\s_n(\setT)} :=\{w\in\s_n\mid \setT(w)=\setT\}$.
Let \Blue{$\m_n(\setT)$} be those bi-leveled trees whose order
ideal consists of the nodes in $\setT$. 
Note that $\beta(\s_n(\setT))=\m_n(\setT)$ and $\beta^{-1}(\m_n(\setT))=\s_n(\setT)$.

\begin{lemma}\label{lem:iota}
 The map $\iota \colon\ \m_n(\setT) \to \s_n(\setT)$ is a map of posets.
\end{lemma}

\begin{proof}
 Let $\setT=\{ 1=T_1<\dotsb<T_r\}$.
 Setting $T_{r+1}=n+1$, define $a_i:=T_{i+1}-T_i-1$ for $i=1,\dotsc,r$. 
 Then $b\mapsto (t_0,(t_1,\dotsc,t_r))$ gives an isomorphism of posets,
\[
   \m_n(\setT)\ \xrightarrow{\,\sim\,}\ 
   \y_{r-1}\times\y_{a_1}\times\dotsb\times\y_{a_r}\,.
\]
 As the maps $\Min,\Max\colon \y_a\to\s_a$ are order-preserving, the proof of
 Theorem~\ref{thm:beta_preimages} gives the desired result.
\end{proof}

\begin{proof}[Proof of Theorem \ref{thm:iota}]

Let $b\lessdot c$ be a cover in $\m_n$.
We will show that $\iota(b)\leq\iota(c)$ in $\s_n$.
Suppose that $b=(t;\setT)$, with $\setT=\{1=T_1<\dotsb<T_r\}$. 
Let $\iota(b)$ have bi-leveled factorization $\iota(b)=u_1v^1u_2v^2\dotsc u_rv^r$,  
and set $k:=n+1-|\setT|$. 

The result is immediate if the cover $b\lessdot c$ is of type $(i)$, for then $b,c\in\m_n(\setT)$
and $\iota\colon\m_n(\setT)\to\s_n$ is order-preserving, as observed in Lemma~\ref{lem:iota}.

Now suppose that $b\lessdot c$ is a cover of type $(ii)$. 
Set $w:=\iota(b)$. 
We claim that $w\lessdot (k,k{+}1)w$ and $\iota(c)=(k,k{+}1)w$.
Now, $u_1=k$ labels the leftmost node of $b$, so the first claim is immediate.
Note that $u_2$ labels the parent of the node labeled $b$. 
This parent has no other child in $\setT$, so we must have $u_2<u_3$.
As $u_2u_3\dotsc u_r$ is $231$-avoiding and contains $k{+}1$, we must have $u_2=k{+}1$.
This shows that 
\[
   \iota(c) \ =\ (k,k{+}1)w\ =\ u_2\,(\Blue{v^1u_1v^2})\, u_3\dotsc u_r v^r\,.
\]
Indeed,  $u_2$ is minimal among $u_2,\dotsc,u_r$ and $u_3\dotsc u_r$ is 
$231$-avoiding, thus $\min_0(c)=u_2\dotsc u_r$.
The bi-leveled factorization of $(k,k{+}1)w$ gives 
$(v^1u_1v^2, v^3,\dotsc,v^r)$, which we claim is $\underline{\Max}(c)$.
As $u_1$ is the largest letter in the sequence, we need only check that
$v^1u_1v^2$ is $132$-avoiding.
But this is true for $v^1$ and $v^2$ and there can be no $132$-pattern involving
$u_1$ as the letters in $v^1$ are all greater than those in $v^2$.

Finally, suppose that $b\lessdot c$ is of type $(iii)$. 
Then $c=(t;\setT\setminus\{T_j\})$ for some $j>2$.
We will find a permutation $w'\in\beta^{-1}(b)$ satisfying $(k,k{+}1)w'\in\beta^{-1}(c)$ and 
\begin{equation}\label{Eq:chain}
   \iota(b)\ \leq\ w'\ \lessdot (k,k{+}1)w'\ \leq\ \iota(c)\,.
\end{equation}
Let $w'\in\beta^{-1}(b)$ be the minimal permutation having bi-leveled factorization 
\[
  w'\ =\ u'_1 v^1 u'_2  \dotsc u'_r v^r\,, \quad\mbox{with}\quad u'_j=k{+}1\,.
\]
Here $(v^1,\dotsc,v^r)=\underline{\Max}(b)$ is the same sequence as in $\iota(b)$. 
The structure of $\beta^{-1}(b)$ implies that $\iota(b)\leq w'$.
We also have
\[
  w'\ \lessdot\ (k,k{+}1)w'
  \qquad\mbox{and}\qquad
  \beta((k,k{+}1)w')\ =\ c\,.
\]
While $\iota(c)$ and $(k,k{+}1)w'$ are not necessarily equal, we do have that
\[
   (k,k{+}1)w'|_{\setT\setminus\{T_j\}}\ =\ 
    u'_{j}u'_2\dotsc u'_{j-1} u'_{j+1} \dotsc u'_r
\]
and $u'_2\dotsc u'_{j-1} u'_{j+1}\dotsc u'_r$ is $231$-avoiding. 
That is, $(k,k{+}1)w'|_{\setT\setminus\{T_j\}} = \iota(c)|_{\setT\setminus\{T_j\}}$. 
Otherwise, $w'$ would not be minimal. 
The bi-leveled factorization of $(k,k{+}1)w'$ is 
\[
  u'_j\, v^1\, u'_2\, \dotsc\, u'_{j-1}\, (v^{j-1}u'_1v^{j})\, 
  u'_{j+1}\,\dotsc\, u'_r \, v^r\,,
\]
and we necessarily have 
$(v^1,\dotsc,v^{j-1}u'_1v^{j},\dotsc,v^r)\leq\underline{\Max}(c)$, which imples
that $ (k,k{+}1)w'\leq \iota(c)$. 
We thus have the chain \eqref{Eq:chain} in $\s_n$,
completing the proof.
\end{proof}

If $b\lessdot c$ is the cover of type $(iii)$ in Figure~\ref{Fig:M-cover}, the
chain~\eqref{Eq:chain} from $\iota(b)$ to $\iota(c)$ is
\[
   \Blue{4}3\Blue{5}\Blue{7}12\Blue{6}\ \leq\ 
   \Blue{4}3\Blue{6}\Blue{7}12\Blue{5}\ \lessdot\ 
   \Blue{5}3\Blue{6}\Blue{7}124\ \leq\ 
   \Blue{5}4\Blue{6}\Blue{7}123\,.
\]

\subsection{Tree enumeration}\label{sec: enumeration}

Let
 \begin{align}
\notag
  \Blue{\enum{q}S} \ &:=\ \sum_{n\geq 0} n! q^n\ =\ 1 + q + 2q^2 + 6q^3 + 24q^4 + 120q^5+\dotsb 
\intertext{be the enumerating series of permutations, and define $\enum{q}M$ and
$\enum{q}Y$  similarly}
%
  \Blue{\enum{q}M} \ &:=\ \sum_{n\geq0} A_n q^n \ =\ 1 + q + 2q^2 + 6q^3 + 21q^4 + 80q^5 + \dotsb \
  , \\
\notag
  \Blue{\enum{q}Y} \ &:=\ \sum_{n\geq0} C_n q^n \ =\ 1 + q + 2q^2 + 5q^3 + 14q^4 + 42q^5 + \dotsb\ , 
\end{align}
where $\Blue{A_n}:=|\m_n|$ and $C_n:=|\y_n|$ 
are the Catalan numbers $\frac{1}{n+1}\binom{2n}{n}$, whose enumerating series satisfies
\[
	\enum{q}Y\ =\ \frac{1-\sqrt{1-4q}}{2q}\ =\ \frac{2}{1+\sqrt{1-4q}}\ .
\]
Bi-leveled trees are Catalan-like~\cite[Theorem 3.1]{For:2008}: 
for $n\geq 1$, $A_n = C_{n-1} + \sum_{i=1}^{n-1} A_i\,A_{n-i}$.
See also \cite[A121988]{Slo:oeis}. 
Their enumerating series satisfies
\[
	\Blue{\enum{q}M}\ = \
		1 + q \enum{q}Y \cdot \enum{q\enum{q}Y}Y .
\]
We will also be interested in 
$\Blue{\enum{q}{M_+}} := \sum_{n>0}A_n q^n = q \enum{q}Y \cdot \enum{q\enum{q}Y}Y$.

\begin{theorem}\label{thm:hilbert_quotients} 
   The only nontrivial quotients of the enumerating series
   $\enum{q}S$, $\enum{q}{M}$, $\enum{q}{M_+}$, and $\enum{q}{Y}$
   whose expansions have nonnegative coefficients are
\[
    \enum{q}S/\enum{q}M, \quad \enum{q}S/\enum{q}Y, \quad \enum{q}{M_+}/\enum{q}Y, 
    \quad\hbox{and}\quad \enum{q}M/\enum{q}Y .
\]
\end{theorem}

\begin{proof}
 We prove the positivity of the quotient $\enum{q}S/\enum{q}M$ in
 Section~\ref{sec:Hopf_module_msym}. 
 The positivity of $\enum{q}S/\enum{q}Y$ was established
 after~\cite[Theorem~7.2]{AguSot:2006}, which shows that $\ssym$ is a smash product over $\ysym$.  

 For the positivity of $\enum{q}{M_+}/\enum{q}Y$, we use~\cite[Proposition 3]{Bar:2005}, which 
 computes $\enum{q\enum{q}Y}Y = \sum_{n>0} B_n q^{n-1}$, where 
 \begin{equation}\label{eq: Bn}
   B_1\ :=\ C_0 \quad \hbox{and} \quad
   B_n\ :=\ \sum_{k=0}^{n-1} \frac{k}{n-1}\binom{2n-k-3}{n-k-1} C_k \quad \hbox{for }  n>1\,.
 \end{equation}
 In particular, $B_n\geq0$ for all $n\geq0$. Returning to the quotient, we have
 \begin{equation*}
   \frac{\enum{q}{M_+}}{\enum{q}Y}\ =\
   \frac{q \enum{q}Y \cdot \enum{q\enum{q}Y}Y}{\enum{q}Y}\ =\ q\enum{q\enum{q}Y}Y\,,
 \end{equation*}
 so $\enum{q}{M_+}/ \enum{q}Y=\sum_{n>0}B_nq^n$ has nonnegative coefficients. 

 For $\enum{q}M /\enum{q}Y$, use the identity $1/\enum{q}Y = 1 - q\enum{q}Y$ to obtain
 \[
    \frac{\enum{q}M}{\enum{q}Y}\ =\ \enum{q}{M_+} + 1 - q\enum{q}Y \ =\ 
    1 + \sum_{n>0} (B_n - C_{n-1})q^n\,.
 \]
 Positivity is immediate as $B_n-C_{n-1} \geq0$ for $n>0$.

 We leave the proof that the remaining quotients have negative coefficients to the reader's
 computer. 
\end{proof}

\begin{remark}
 Up to an index shift, the quotient $\enum{q}{M_+}/\enum{q}Y$ corresponds to the
 sequence~\cite[A127632]{Slo:oeis} begining  with $(1,1,3,11,44,185,804)$. 
 We give a new combinatorial interpretation of this sequence in Corollary~\ref{cor:free}.
\end{remark}

\section{The Algebra $\msym$}\label{sec: msym}

Let $\Blue{\msym}:=\bigoplus_{n\geq0} {\msym}_{n}$ denote the graded $\Q$--vector space whose $n^{th}$
graded piece has the basis $\{F_b \mid b\in \m_n\}$.  
The maps $\beta\colon\s_\bb\to\m_\bb$ and $\phi\colon\m_\bb\to\y_\bb$ of graded sets
induce surjective maps of graded vector spaces
 \begin{equation}\label{Eq:tau_factor}
   \ssym\ \xrightarrow{\ \Blue{\bbeta}\ }\ 
   \msym\ \xrightarrow{\ \Blue{\bphi}\ }\ \ysym
    \qquad F_w\ \mapsto\ F_{\beta(w)}\ \mapsto\ F_{\phi(\beta(w))}\,,
 \end{equation}
which factor the Hopf algebra map $\btau\colon\ssym\to\ysym$,
as $\phi(\beta(w))=\tau(w)$.
We will show how the maps $\bbeta$ and $\btau$ induce on $\msym$ the structures of an 
algebra, of a $\ssym$-module, and of a $\ysym$-comodule so that the
composition~\eqref{Eq:tau_factor} factors the map $\btau$ as maps of algebras, of
$\ssym$-modules, and of $\ysym$-comodules.
%
%

\subsection{Algebra structure on $\msym$}\label{sec:algebra}
For $b,c\in\m_\bb$ define
 \begin{equation}\label{Eq:multiplication}
   F_b \bm{\cdot} F_c\ =\ \bbeta( F_w \bm{\cdot} F_v)\,,
 \end{equation}
where $w,v$ are permutations in $\s_\bb$ with $b=\beta(w)$ and $c=\beta(v)$.

\begin{theorem}\label{thm:Msym_algebra}
  The operation $F_b \bm{\cdot} F_c$ defined by~\eqref{Eq:multiplication} is independent
  of choices of $w,v$ with $\beta(w)=b$ and $\beta(v)=c$ and it endows $\msym$ with the
  structure of a graded connected algebra such that the map $\bbeta\colon\ssym\to\msym$ is
  a surjective map of graded connected algebras.
\end{theorem}

If the expression $\bbeta(F_w\cdot F_v)$ is independent of choice of $w\in\beta^{-1}(b)$ 
and $v\in\beta^{-1}(c)$, then the map $\bbeta$ is automatically multiplicative. 
Associative and unital properties for $\msym$ are then inherited from those for $\ssym$, 
and the theorem follows. 
To prove independence (in Lemma \ref{thm: independence}), we formulate a description of \eqref{Eq:multiplication} in terms of splittings and graftings of bi-leveled trees. 

Let $s\psplit(s_0,\dotsc,s_m)$ be a splitting on the underlying tree of a bi-leveled
tree $b=(s;\setS)\in\m_n$. 
Then the nodes of $s$ are distributed among the nodes of the partially ordered forest
$(s_0,\dotsc,s_m)$ so that the order ideal $\setS$ gives a sequence of 
order ideals in the trees $s_i$. 
Write $b\psplit(b_0,\dotsc,b_m)$ for the corresponding splitting of the bi-leveled tree $b$, 
viewing $b_i$ as $(s_i;\setS|_{s_i})$. 
(Note that only $b_0$ is guaranteed to be a bi-leveled tree.) 
Given $c=(t;\setT)\in\m_m$ and a splitting $b\psplit(b_0,\dotsc,b_m)$ of $b\in\m_n$, form
a bi-leveled tree \Blue{$(b_0,\dotsc,b_m)/c$} whose underlying tree is $(s_0,\dotsc,s_m)/t$
and whose order ideal is either
 \begin{equation}\label{Eq:split_graft}
  \hskip1.3cm\begin{minipage}[c]{10cm}
   \begin{enumerate}
    \item[$({i)}$\ ] $\setT$, \ if $b_0\in \m_0$, or

    \item[$({ii)}$\ ]  $\setS\cup\{\mbox{the nodes of $t$}\}$, \ if $b_0\not\in \m_0$.
   \end{enumerate}
  \end{minipage}
 \end{equation}
%
%
\begin{lemma}\label{thm: independence}
  The product~$\eqref{Eq:multiplication}$  is independent of choices of $w,v$ 
  with $\beta(w)=b$ and $\beta(v)=c$. 
  For $b\in\m_n$ and $c\in\m_m$, we have
\[
     F_b \bm{\cdot} F_c\ =\ 
     \sum_{b\psplit(b_0,\dotsc,b_m)} F_{(b_0,\dotsc,b_m)/c} \,.
\]
\end{lemma}

\begin{proof}
Fix any $w\in\beta^{-1}(b)$ and $v\in\beta^{-1}(c)$.
The bi-leveled tree $\beta((w_0,\dotsc,w_m)/v)$ associated to a splitting
$w\psplit(w_0,\dotsc,w_m)$ has underlying tree $(s_0,\dotsc,s_m)/t$, where
$s\psplit(s_0,\dotsc,s_m)$ is the 
induced splitting on the underlying tree $s=\tau(w)=\phi(b)$.
Each node of $(w_0,\dotsc,w_m)/v$ comes from a node of either $w$ or $v$, with the
labels of nodes from $w$ all smaller than the labels of nodes from $v$.
Consequently, the leftmost node of $(w_0,\dotsc,w_m)/v$ comes from either 
 \begin{enumerate}
  \item[$({i)}$] $v$, and then $\setT((w_0,\dotsc,w_m)/v)=\setT(v)=\setT(c)$, or

  \item[$({ii)}$] $w$, and then 
         $\setT((w_0,\dotsc,w_m)/v)=\setT(w)=\setT(b)\cup\{\mbox{the nodes of $v$}\}$.
 \end{enumerate}
The first case is when $w_0\in\s_0$ and the second case is when $w_0\not\in\s_0$. 
\end{proof}

Here is the product 
$F_{\includegraphics{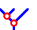}}
  \bm\cdot F_{\includegraphics{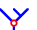}}$, 
together with the corresponding splittings of $\includegraphics{figures/m12.s.eps}$,
\newcommand{\SPA}{\hspace{1.665em}}
\newcommand{\spa}{\hspace{2.24em}}
 \begin{align*}
 F_{\includegraphics{figures/m12.s.eps}} 
 \bm\cdot F_{\includegraphics{figures/m21.s.eps}} 
 \ &=\
     {F_{\includegraphics{figures/m1243.s.eps}}} + 
     {F_{\includegraphics{figures/m1423.s.eps}}} + 
     {F_{\includegraphics{figures/m1432.s.eps}}} +
     {F_{\includegraphics{figures/m4123.s.eps}}} + 
     {F_{\includegraphics{figures/m4231.s.eps}}} + 
     {F_{\includegraphics{figures/m4312.s.eps}}}\,.\\
 &\spa\includegraphics{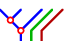}\rule{0pt}{16pt}
   \SPA\includegraphics{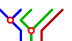}
   \SPA\includegraphics{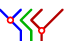}
   \SPA\includegraphics{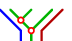}
   \SPA\includegraphics{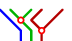}
   \SPA\includegraphics{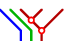}
 \end{align*}
%

\subsection{$\ssym$ module structure on $\msym$}\label{sec: ssym-mod} 
Since $\bbeta$ is a surjective algebra map, $\msym$ becomes a $\ssym$-bimodule with the
action
\[
   F_w \bm\cdot F_b \bm\cdot F_v\ =\  F_{\beta(w)}\bm\cdot F_b\bm\cdot F_{\beta(v)}\ .
\]
The map $\btau$ likewise induces on $\ysym$ the structure of a $\ssym$-bimodule, and the
maps $\bbeta$, $\bphi$, and $\btau$ are maps of $\ssym$-bimodules.

Curiously, we may use the map $\iota\colon\m_\bb\to\s_\bb$ to define the structure of a right
$\ssym$-comodule on $\msym$, 
\[
   F_b\ \longmapsto\ \sum_{\iota(b)\psplit(w_0,w_1)} F_{\beta(w_0)}\otimes F_{w_1}\,.
\]
This induces a right comodule structure, because if
$\iota(b)\psplit(w_0,w_1)$, then $w_0=\iota(\beta(w_0))$, which may be checked using the 
characterization of $\iota$ in terms of pattern avoidance, as explained in Remark~\ref{Rem:pattern}.

While $\msym$ is both a right $\ssym$-module and right $\ssym$-comodule, it is not an 
$\ssym$--Hopf module. 
For if it were a Hopf module, then the fundamental theorem of Hopf modules (see
Remark~\ref{Rem:covariants}) would imply that the series $\enum{q}{M}/\enum{q}{S}$ has positive 
coefficients, which contradicts Theorem~\ref{thm:hilbert_quotients}. 

\subsection{$\ysym$-comodule structure on $\msym$}\label{sec: ysym-mod} 
For $b\in\m_\bb$, define the linear map $\brho\colon\msym \to \msym\otimes \ysym$ by
 \begin{equation}\label{Eq:coaction}
    \brho (F_b)\ =\ \sum_{b \psplit (b_0,b_1)} F_{b_0}\otimes F_{\phi(b_1)}\,.
 \end{equation}
By $\phi(b_1)$, we mean the tree underlying $b_1$.

\begin{example}\label{ex: ysym coaction}
In the fundamental bases of $\msym$ and $\ysym$, we have
\[
   \brho(F_{\includegraphics{figures/m2431.s.eps}})\  =\  
        {F_{\includegraphics{figures/m2431.s.eps}} \otimes 1} 
      + {F_{\includegraphics{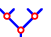}} \otimes 
         F_{\includegraphics{figures/1.eps}}}
      + {F_{\includegraphics{figures/m12.s.eps}} \otimes 
         F_{\includegraphics{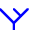}}} 
      + {F_{\includegraphics{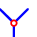}} \otimes
         F_{\includegraphics{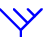}}}
      + {1 \otimes F_{\includegraphics{figures/3421.s.eps}}}\, .
\]
\end{example}

\begin{theorem}\label{thm:Ysym_Coaction}
 Under $\brho$, $\msym$ is a right $\ysym$-comodule.
\end{theorem}

\begin{proof} 
 This is counital as $(b,\includegraphics{figures/0.eps})$ is a splitting of $b$.
 Coassociativity is also clear as both $(\brho\otimes 1)\brho$ and $(1\otimes \Delta)\brho$
 applied to  $F_b$ for $b\in\m_\bb$ are sums of  terms 
 $F_{b_0}\otimes F_{\phi(b_1)}\otimes F_{\phi(b_2)}$
 over all splittings $b\psplit(b_0,b_1,b_2)$.
\end{proof}

Careful bookkeeping of the terms in $\brho(F_b \cdot F_c)$ show that it equals 
$\brho(F_b)\cdot \brho(F_c)$ for all $b,c\in\m_\bb$ and thus $\msym$ is a
$\ysym$--comodule algebra.  
Hence, $\bphi$ is a map of $\ysym$--comodule algebras,
and in fact $\bbeta$ is also a map of $\ysym$--comodule algebras.
We leave this to the reader, and will not pursue it further.

Since $\btau\colon\ssym\to\ysym$ is a map of Hopf algebras, $\ssym$ is naturally a right
$\ysym$-comodule where the comodule map is the composition
\[
  \ssym\ \xrightarrow{\ \Delta\ }\ 
  \ssym \otimes\ssym\ \xrightarrow{\ 1\otimes\btau\ }\ 
  \ssym \otimes\ysym\,.
\]
With these definitions, the following lemma is immedate.

\begin{lemma}\label{lemma:comodule_map}
 The maps $\btau$ and $\bphi$ are maps of right $\ysym$-comodules.
\end{lemma}

In particular, we have the equality of maps $\ssym \to \msym \otimes\ysym$,
 \begin{equation}\label{Eq:comod_map}
   \brho\circ\bbeta\ =\ (\bbeta\otimes\btau)\circ \Delta\,.
 \end{equation}
%

\subsection{Coaction in the monomial basis}\label{sec: coaction in M basis}

The coalgebra structures of $\ssym$ and $\ysym$ were elicudated by considering a second basis
related to the fundamental basis via M\"obius inversion.
For $b\in \m_n$, define 
 \begin{equation}\label{eq: msym M-basis}
	\Blue{M_b}\ :=\ \sum_{b\leq c} \mu(b,c) F_{c}\,,
 \end{equation}
where $\mu(\cdot,\cdot)$ is the M\"obius function on the poset $\m_n$. 

Given $b\in \m_m$ and  $s\in \y_q$, write
\Blue{$b\backslash s$} for the bi-leveled tree with $p+q$ nodes whose underlying tree is
formed by grafting the root of $s$ onto the rightmost leaf of $b$, but whose order ideal is
that of $b$.  
Here is an example of $b$, $s$, and $b\backslash s$,
\[
   \includegraphics{figures/m2143.d.eps} \qquad \ 
   \includegraphics{figures/21.d.eps} \qquad \ 
   \includegraphics{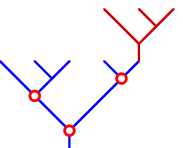} 
   \ \ \raisebox{5pt}{$=$}\ \
   \includegraphics{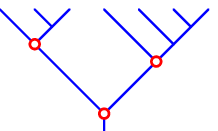}
   \raisebox{1pt}{.}
\]
Observe that we cannot have $b=\includegraphics{figures/0.eps}$ in this construction.

The maximum bi-leveled tree with a given underlying tree $t$
is $\beta(\Max(t))$,  which has order ideal $\setT$ consisting only of the nodes of
$t$ along its leftmost branch. 
Here are three such trees of the form $\beta(\Max(t))$,
\[
  \includegraphics{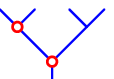} \qquad
  \includegraphics{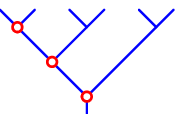} \qquad
  \includegraphics{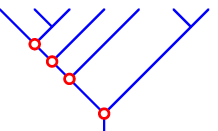} 
   \raisebox{1pt}{.}
\]

\begin{theorem}\label{thm: M-basis coaction} 
 Given $b=(t;\setT)\in \m_\bb$, we have 
 \[ 
  \brho(M_b)\ =\ \left\{\begin{array}{lcl}
    {\displaystyle \sum_{b=c\backslash s}M_c\otimes M_s} 
        &\ &\mbox{if }b\neq\beta(\Max(t))\\
    {\displaystyle \sum_{b=c\backslash s} M_c \otimes M_s} 
     \ +\  1 \otimes M_t &\ &\mbox{if }b=\beta(\Max(t))\rule{0pt}{20pt}
  \end{array}\right.\ .
 \] 
\end{theorem}

For example, 
 \begin{align*}
\brho(M_{\includegraphics{figures/m1432.s.eps}}) \ &=\ 
    {M_{\includegraphics{figures/m1432.s.eps}} \otimes 1}  \\
\brho(M_{\includegraphics{figures/m2431.s.eps}}) \ &=\ 
    {M_{\includegraphics{figures/m2431.s.eps}} \otimes 1} 
	 + {M_{\includegraphics{figures//m132.s.eps}} \otimes 
	M_{\includegraphics{figures/1.eps}}} \\
\brho(M_{\includegraphics{figures/m3421.s.eps}}) \ &=\ 
     {M_{\includegraphics{figures/m3421.s.eps}} \otimes 1} 
	 + {M_{\includegraphics{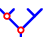}} \otimes 
	M_{\includegraphics{figures/1.eps}}}
	 + {M_{\includegraphics{figures/m12.s.eps}} \otimes 
	M_{\includegraphics{figures/21.s.eps}}} 
	 \ +\ {1 \otimes M_{\includegraphics{figures/3421.s.eps}}} .
\end{align*}

Our proof of Theorem~\ref{thm: M-basis coaction} 
uses Proposition~\ref{prop:M-basis} and the following results.

\begin{lemma}\label{lemma:Beta_M}
 For any bi-leveled tree $b\in \m_\bb$, we have
\[
    \bbeta\biggl(\;\sum_{\beta(w)=b} M_w \biggr) = M_b \,.
\]
\end{lemma}

\begin{proof}
 Expand the left hand side in terms of the fundamental bases to get
\[
   \bbeta\biggl( \sum_{\beta(w)=b}\ \sum_{w\leq v}  \mu_\s(w,v) F_v \biggr)
   \ =\   
   \sum_{\beta(w)=b}\ \sum_{w \leq v} \mu_\s(w,v) F_{\beta(v)} \,.
\]
 As $\beta$ is surjective, we may change the index of summation to $b\leq c$ in $\m_\bb$ to
 obtain
\[
   \sum_{b\leq c} \biggl(\;\sum_{\substack{\beta(w)=b\\\beta(v)=c}} \mu_{\s}(w,v)\biggr) F_{c}\,.
\]
 By Theorems~\ref{thm: subGalois} and~\ref{thm:iota},
 the inner sum is $\mu_{\m}(b,c)$, so this sum is $M_b$.
\end{proof}

Recall that $w=u\backslash v$ only if $\tau(w)=\tau(u)\backslash\tau(v)$
and the values of $w$ in the nodes of $u$ exceed the values in the nodes of $v$.
We always have the trivial decomposition $w=(\emptyset,w)$.
Suppose that $w=u\backslash v$ with $u\neq\emptyset$ a nontrivial decomposition. 
If $\beta(w)=b=(t;\setT)$, then $\setT$ is a subset of the nodes of $u$ so that 
$\beta(u)=(\tau(u);\setS)$ and $b=\beta(u)\backslash \tau(v)$.
Moreover, for every decomposition $b=c\backslash s$ and every $u,v$ with 
$\beta(u)=c$ and $\tau(v)=s$, we have $b=\beta(u\backslash v)$.
Thus, for $b\in\m_\bb$, we have
 \begin{equation}\label{Eq:disjoint_decomposition}
    \bigsqcup_{\beta(w)=b}\ 
    \bigsqcup_{\substack{w=u\backslash v\\u\neq\emptyset}}\; (u,v)
   \ =\ 
    \bigsqcup_{b=c\backslash t}\ 
    \bigsqcup_{\beta(u)=c}\ 
    \bigsqcup_{\tau(v)=t}\; (u,v)\ .
 \end{equation}

\begin{proof}[Proof of Theorem$~\ref{thm: M-basis coaction}$]
 Let $b=(t;\setT)$ with $t\neq \includegraphics{figures/0.eps}$.
 Using Lemma~\ref{lemma:Beta_M}, we have
\[
  \brho(M_b)\ =\ 
  \brho\bbeta\biggl(\sum_{\beta(w)=b} M_w \biggr)\ =\ 
  \sum_{\beta(w)=b}\brho\bbeta M_w\ .
\]
By~\eqref{Eq:comod_map}, \eqref{Eq:disjoint_decomposition}, 
and~\eqref{Eq:SSym_M_coprod}, this equals 
 \begin{multline*}
   \sum_{\beta(w)=b}\ \sum_{\substack{w=u\backslash v\\u\neq\emptyset}}
     \bbeta(M_u)\otimes \btau(M_v)\quad +\ 
    \sum_{\beta(w)=b} \bbeta(M_\emptyset)\otimes \btau(M_w)\\
  =\ \sum_{b=c\backslash s} 
   \biggl(\;\sum_{\beta(u)=c} \bbeta(M_u)\biggr) \otimes 
   \biggl(\;\sum_{\tau(v)=s} \btau(M_t)\biggr)\quad +\ 
    \sum_{\beta(w)=b} 1\otimes \btau(M_w)\ .
 \end{multline*}
By Lemma~\ref{lemma:Beta_M} and~\eqref{eq:tau_M-basis}, the first sum becomes 
$\sum_{b=c\backslash s}M_c\otimes M_s$ and the second sum vanishes unless $b=\beta(\Max(t))$. 
This completes the proof.
\end{proof}

\section{Hopf Variations}\label{sec: variations}

\subsection{The $\ysym$--Hopf module $\msym_+$}\label{sec:Hopf_module_msym+} 
Let $\Blue{\m_+}:=(\m_n)_{n\geq 1}$ be the bi-leveled trees with at least one internal node
and define $\Blue{\msym_+}$ to be the positively graded part of $\msym$, which has bases
indexed by $\m_+$. 
A \demph{restricted splitting} of $b\in\m_+$ is a splitting 
$b\rsplit(b_0,\ldots,b_m)$ with $b_0\in\m_+$, i.e., $b_0\neq \includegraphics{figures/0.eps}$.
Given $b\rsplit(b_0,\ldots,b_m)$ and $t\in\y_m$, form the bi-leveled tree
$(b_0,\dotsc,b_m)/t$ by grafting the ordered forest $(b_0,\dotsc,b_m)$ onto the leaves of
$t$, with order ideal consisting of the nodes of $t$ together with the nodes of the forest
coming from the order ideal of $b$, as in~\eqref{Eq:split_graft}({\it ii}).

We define an action and coaction of $\ysym$ on $\msym_+$ that are similar to the 
product and coaction on $\msym$. 
They come from a second collection of polytope maps $\m_n \onto \y_{n-1}$ arising from viewing the vertices of $\m_n$ as \emph{painted trees} on $n-1$ nodes (see \cite{BoaVog:1973,For:2008}). 
For $b\in \m_{+}$ and $t\in \y_m$, set 
 \begin{equation}\label{eq: msym+ action and coaction}
  \begin{array}{lcl}
   F_b\bm\cdot F_t  &=&
     {\displaystyle \sum_{b \rsplit (b_0,\dotsc,b_m)} F_{(b_0,\dotsc,b_m)/t}}\,,\\
    \Blue{\brho_+} (F_b) &=&
      {\displaystyle \sum_{b \rsplit (b_0,b_1)} F_{b_0}\otimes F_{\phi(b_1)}}\,.\rule{0pt}{20pt}
  \end{array}
 \end{equation}

For example, in the fundamental bases of $\msym_+$ and $\ysym$, we have
\begin{align*}
  F_{\includegraphics{figures/m21.s.eps}}  
  \bm\cdot F_{\includegraphics{figures/21.s.eps}} \ &=\ 
  {F_{\includegraphics{figures/m2143.s.eps}}} 
+ {F_{\includegraphics{figures/m2413.s.eps}}} 
+ {F_{\includegraphics{figures/m2431.s.eps}}}\,,\\
  \brho_+(F_{\includegraphics{figures/m3241.s.eps}}) \ &=\ 
     {F_{\includegraphics{figures/m3241.s.eps}} \otimes 1} + 
     {F_{\includegraphics{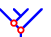}} \otimes  
       F_{\includegraphics{figures/1.eps}}} +   
     {F_{\includegraphics{figures/m21.s.eps}} \otimes 
        F_{\includegraphics{figures/21.s.eps}}} + 
     {F_{\includegraphics{figures/m1.s.eps}} \otimes 
        F_{\includegraphics{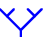}}}\, . 
 \end{align*}

\begin{theorem}\label{thm:msym_+}
The operations in \eqref{eq: msym+ action and coaction} define a $\ysym$--Hopf module structure on $\msym_+$.
\end{theorem}

\begin{proof}
The unital and counital properties are immediate. 
We check only that the action is associative, the coaction is coassociative, and the two structures commute with each other.

{\em Associativity.\,}~%
Fix $b=(t;\setT) \in \m_+$, $r\in \y_m$, and $s\in \y_n$. 
A term in the expression $(F_b\cdot F_r)\cdot F_s$ corresponds to a restricted splitting and grafting 
$b\rsplit (b_0,\dotsc,b_m) \leadsto (b_0,\dotsc, b_m)/r=c$, followed by another 
$c\rsplit (c_0,\dotsc,c_n) \leadsto (c_0,\dotsc,c_n)/t$. 
The order ideal for this term equals $\setT \cup \{\mbox{the nodes of $r$ and $s$}\}$. 
Note that restricted splittings of $c$ are in bijection with pairs of splittings 
\[
	\Bigl(b \rsplit (b_0,\dotsc,b_{m{+}n})\,,\, r \psplit (r_0,\dotsc,r_n)\Bigr).
\] 
Terms of $F_b \cdot (F_r\cdot F_s)$ also correspond to these pairs of splittings. 
The order ideal for this term is again $\setT \cup \{\mbox{the nodes of $r$ and $s$}\}$. That is, 
$(F_b\cdot F_r)\cdot F_s$ and $F_b\cdot (F_r \cdot F_s)$ agree term by term.

{\em Coassociativity.\,}~%
Fix $b=(t;\setT) \in \m_+$. Terms $F_c\otimes F_r\otimes F_s$ in $(\brho_+ \otimes \id)\brho_+(F_b)$ and $(\id \otimes \Delta)\brho_+(F_b)$ both correspond to restricted splittings $b\rsplit(c,c_1,c_2)$, where $\phi(c_1)=r$ and $\phi(c_2)=s$. 
In either case, the order ideal on $c$ is $\setT|_c$. 

{\em Commuting structures.\,}~%
Fix $b=(s;\setS) \in \m_+$ and $t\in \y_m$. 
A term $F_{c_0}\otimes F_{\phi(c_1)}$ in $\brho_+(F_b \cdot F_t)$ corresponds to a choice of 
a restricted splitting and grafting $b \rsplit (b_0,\dotsc,b_m) \leadsto (b_0,\ldots, b_r)/t=c$, 
followed by a restricted splitting $c \rsplit (c_0,c_1)$. 
The order ideal on $c_0$ equals the nodes of $c_0$ inherited from $\setS$, 
together with the nodes of $c_0$ inherited from $t$. 
The restricted splittings of $c$ are in bijection with pairs of splittings 
$\bigl(b \rsplit (b_0,\dotsc,b_{m{+}1}), t \psplit (t_0,t_1)\bigr)$. 
If $t_0\in \y_n$, then the pair of graftings $c_0=(b_0,\dotsc,b_n)/t_0$ and $c_1=(b_{n{+}1},\dotsc,b_m)/t_1$ are precisely the terms appearing in $\brho_+(F_b)\cdot \Delta(F_t)$. 
\end{proof}

The similarity of \eqref{eq: msym+ action and coaction} to the coaction \eqref{Eq:coaction} of $\ysym$ on $\msym$ gives the following result, whose proof we leave to the reader.

\begin{corollary}\label{cor:brhoPlus}
 For $b\in\m_+$, we have
\[ 
   \brho_+(M_b)\ =\ \sum_{{b=c\backslash s}} M_{c} \otimes M_s\,.
\] 
\end{corollary}

This elucidates the structure of $\msym_+$.
Let $\Blue{\calB}\subset \m_+$ be the indecomposable bi-leveled trees---those with 
only trivial decompositions, $b=b\backslash \includegraphics{figures/0.eps}$\,.
Then $(t;\setT)\in\calB$ if and only if $\setT$ contains the rightmost node of $t$.
Every tree $c$ in $\m_+$ has a unique decomposition  $c=b\backslash s$ where $b\in\calB$ and
$s\in\y_\bb$. 
Indeed, pruning $c$ immediately above the rightmost node in its order ideal
gives a decomposition $c=b\backslash s$ where $b\in\calB$ and
$s\in\y_\bb$. 
This induces a bijection of graded sets,
\[
   \m_+\ \longleftrightarrow\ \calB\times\y_\bb\,.
\]
Moreover, if $b\in\calB$ and $s\in\y_\bb$, then Corollary~\ref{cor:brhoPlus}
and~\eqref{Eq:YSym_M_coprod} together imply that
 \begin{equation}\label{Eq:new_brhoPlus}
   \brho_+(M_{b\backslash s})\ =\ 
   \sum_{s=r\backslash t} M_{b\backslash r}\otimes M_t\,.
 \end{equation}

Note that $\Q\calB \otimes\ysym$ is a graded right $\ysym$-comodule
with structure map,
\[
    b \otimes M_s\ \longmapsto\ b\otimes (\Delta M_s)\,,
\]
for $b\in\calB$ and $s\in\y_\bb$.
Comparing this with~\eqref{Eq:new_brhoPlus}, we deduce the following algebraic and combinatorial
facts.

\begin{corollary}\label{cor:free}
 The map $\Q\calB\otimes\ysym\to\msym_+$ defined by 
 $b\otimes M_s\mapsto M_{b\backslash s}$ is an isomorphism of graded right $\ysym$ comodules.

 The quotient of enumerating series $\enum{q}M_+/\enum{q}Y$ is equal to the enumerating series
 of the graded set $\calB$.
\end{corollary}

In particular, if $\calB_n:=\calB\cap \m_n$, then  
$|\calB_n|=B_n$ by~\eqref{eq: Bn}.

\begin{remark}\label{Rem:covariants}
 The \demph{coinvariants} in a right comodule $M$ over a coalgebra $C$ are 
 $\Blue{M^{\mathrm{co}}} := \left\{ m \in M \mid \brho(m) = m\otimes 1 \right\}$.
 We identify the vector space $\Q\calB$ with $\msym_+^{\mathrm{co}}$ via
 $b\mapsto M_b$.
 The isomorphism $\Q\calB\otimes\ysym\to\msym_+$ is a special case of the
 Fundamental Theorem of Hopf Modules~\cite[Theorem~1.9.4]{Mont:1993}:
 If $M$ is a Hopf module over a Hopf algebra $H$, then 
 $M \simeq M^{\mathrm{co}} \otimes H$ as Hopf modules. 
\end{remark}

%
\subsection{Hopf module structure on $\msym$}
\label{sec:Hopf_module_msym}

We use Theorem~\ref{thm: M-basis coaction}  to identify the $\ysym$-coinvariants in $\msym$.
Let $\Blue{\calB'}$ be those indecomposable bi-leveled trees which 
are not of the form $\beta(\Max(t))$, for some $t\in\y_+$,
together with $\{\includegraphics{figures/0.eps}\}$.

\begin{corollary}\label{cor:coinvariants}
 The $\ysym$-coinvariants of $\msym$ have a basis $\{ M_b\mid b\in \calB'\}$.
\end{corollary} 

For $n>0$, the difference $\calB_n\setminus \calB'_n$ consists of indecomposable bi-leveled
trees with $n$ nodes of the form $\beta(\Max(t))$.
If $\beta(\Max(t))\in\calB_n$, then $t=s\vee\includegraphics{figures/0.eps}$,
for some $s\in\y_{n-1}$, and so $|\calB'_n|=B_n-C_{n-1}$, which we saw in the proof of
Theorem~\ref{thm:hilbert_quotients}.

For $t\in\y_\bb$, set $\includegraphics{figures/0.eps}{\bbslash}t:=\beta(\Max(t))$,
and if $\includegraphics{figures/0.eps}\neq b\in\mathcal B'$, set 
$b{\bbslash}t:=b\backslash t$.
Every bi-leveled tree uniquely decomposes as $b {\bbslash} t$ with $b\in\calB'$
and $t\in\y_\bb$. 
By Theorem~\ref{thm: M-basis coaction}, 
$M_b\otimes M_t\mapsto M_{b{\bbslash} t}$ induces an isomorphism of right $\ysym$-comodules,
 \begin{equation}\label{Eq:trsansfer}
  \msym^{\mathrm{co}}\otimes \ysym\ \longrightarrow\ \msym\,,
 \end{equation}
where the structure map on $\msym^{\mathrm{co}}\otimes \ysym$ is 
$M_b\otimes M_t\mapsto M_b\otimes \Delta(M_t)$.
Treating $\msym^{\mathrm{co}}$ as a trivial $\ysym$-module,
$M_b\cdot M_t=\varepsilon(M_t)M_b$, $\msym^{\mathrm{co}}\otimes \ysym$ becomes a right
$\ysym$-module. 
As explained in~\cite[Example 1.9.3]{Mont:1993}, this makes
$\msym^{\mathrm{co}}\otimes \ysym$ into a $\ysym$--Hopf module.

We express this structure on $\msym$.
Let $b{\bbslash}t\in \m_{\bb}$ and $s\in \y_\bb$, then 
 \begin{equation}\label{eq: msym action and coaction}
   M_{b\bbslash t} \bm\cdot M_s\ =\ 
     \sum_{r \in t\cdot s} M_{b\bbslash r} \,
    \quad\mbox{and}\quad
   \brho\bigl(M_{b\bbslash t}\bigr)\ =\ \sum_{t=r\backslash s} M_{b\bbslash r} \otimes M_{s}\,,
 \end{equation}
where $t\cdot s$ is the set of trees $r$ indexing the product $M_t \cdot M_s$ in $\ysym$.
The coaction is as before, but the product is new.
It is not positive in the fundamental basis,
\[
   F_{\includegraphics{figures/m1.s.eps}} \cdot
   F_{\includegraphics{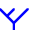}}\ = \	
   F_{\includegraphics{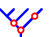}} - 	
   F_{\includegraphics{figures/m132.s.eps}} + 
   F_{\includegraphics{figures/m213.s.eps}} + 
   2F_{\includegraphics{figures/m231.s.eps}}\,.
\]

We complete the proof of  Theorem~\ref{thm:hilbert_quotients}.

\begin{corollary}
   The power series $\enum{q}S / \enum{q}M$ has nonnegative coefficients.
\end{corollary}

\begin{proof}
 Observe that 
\[
  \enum{q}{S} / \enum{q}{M} \ =\  
	  \Bigl(\enum{q}S / \enum{q}Y \Bigr) \big/ 
	 \Bigl(\enum{q}M / \enum{q}Y \Bigr) \,.
\]
 Since both $\ssym$ and $\msym$ are right $\ysym$--Hopf modules, the two quotients of
 enumerating series on the right are generating series for their coinvariants, by the
 Fundamental Theorem of Hopf modules.
 Thus
\[
  \enum{q}{S} / \enum{q}{M} \ =\  
  \enum{q}{S^{\mathrm{co}}} / \enum{q}{M^{\mathrm{co}}}\,,
\]
 where $\enum{q}{S^{\mathrm{co}}}$ and $\enum{q}{M^{\mathrm{co}}}$ are the enumerating series
 for $\ssym^{\mathrm{co}}$ and $\msym^{\mathrm{co}}$. To show that $\enum{q}{S^{\mathrm{co}}} / \enum{q}{M^{\mathrm{co}}}$ is nonnegative, we index bases for these spaces by graded sets $\calS$ and $\calB'$, then establish a bijection $\calB' \times \calS' \to \calS$ for some graded subset $\calS' \subset \calS$.

The set $\calB'$ was identified in Corollary \ref{cor:coinvariants}. 
The coinvariants $\ssym^{\mathrm{co}}$ were given in~\cite[Theorem 7.2]{AguSot:2006}
 as a \emph{left Hopf kernel}. The basis was identified as follows. 
Recall that pemutations $u\in\s_\bb$ may be written uniquely in terms of indecomposables,
 \begin{equation}\label{eq:decomposition}
   u\ =\ u_1\backslash \dotsb \backslash u_r
 \end{equation}
(taking $r=0$ for $u=\emptyset$). Let $\calS\subset\s_\bb$ be those permutations $u$ whose rightmost indecomposable component has
 a $132$-pattern, and thus $u\neq \Max(t)$ for any $t\in\y_+$. 
(Note that $u=\emptyset \in \calS$.) Then $\{M_u\mid u\in\calS\}$ is a basis for $\ssym^{\mathrm{co}}$. 

Fix a section $g\colon\m_\bb\to\s_\bb$ of the map $\beta\colon\s_\bb\to\m_\bb$ and 
define a subset $\calS' \subset \calS$ as follows. 
Given the decomposition $u=u_1\backslash \dotsb \backslash u_r$ in \eqref{eq:decomposition} with $r\geq0$, consider the length $\ell\geq0$ of the maximum initial sequence 
$u_1\backslash \dotsb \backslash u_\ell$ of indecomposables belonging to $g(\calB')$. 
Put $u\in \calS'$ if $\ell$ is even. 
Define the map of graded sets
\[
  \Blue{\kappa}\ \colon\ \calB'\times\calS'\longrightarrow\calS
   \quad\mbox{\ by\ }\quad (b,v)\ \longmapsto\ g(b)\backslash v\,.
\]
 The image of $\kappa$ lies in $\calS$ as the last component of a nontrivial $g(b)\backslash v$ is either $g(b)$ or the last component of $v$, neither of which can be $\Max(t)$ for
 $t\in\y_+$. 

 We claim that $\kappa$ is bijective.
 If $u\in\calS'$, then $u=\kappa(\includegraphics{figures/0.eps},u)$.
 If $u\in\calS\setminus\calS'$,  then $u$ has an odd number of initial components
 from $g(\calB')$.
 Letting its first factor be $g(b)$, we see that $u=g(b)\backslash u'=\kappa(b,u')$ with $u'\in\calS'$.
 This surjective map is injective as the expressions 
 $\kappa(\includegraphics{figures/0.eps},u')$ and 
 $\kappa(b,u')$ with $b\in\calB'_+$ and $u'\in\calS'$ are unique.
 
 This isomorphism of graded sets identifies the enumerating series of the graded set
 $\calS'$ as the quotient $\enum{q}{S^{\mathrm{co}}} / \enum{q}{M^{\mathrm{co}}}$, which
 completes the proof.
\end{proof}

\bibliographystyle{amsplain}
\bibliography{bibl}

\end{document}